\documentclass[
               DIV=15,
               ]{scrartcl}  

\usepackage{color} 
\usepackage{psfrag} 
\usepackage{subcaption}
\usepackage{hyperref}
\usepackage{amsmath}
\usepackage{amssymb}
\usepackage{amsthm}
\usepackage{graphics}
\usepackage{todonotes}
\usepackage{amsfonts}  
\usepackage{adjustbox}    
\usepackage{mathtools}
\usepackage{xcolor}        
\usepackage{blkarray}        
\usepackage{footnote}
\makesavenoteenv{tabular}
\makesavenoteenv{table}
\usepackage{stmaryrd}
\usepackage{bm}
\usepackage{makecell}

\usepackage[square,numbers,comma]{natbib}

\hypersetup{colorlinks,citecolor=blue}

\usepackage{standalone}
 
\usepackage{pgfplots}
\usepackage{pgfplotstable}
\usetikzlibrary{pgfplots.groupplots}
\usetikzlibrary{fit}
\usepackage{filecontents}
\usepackage{tikz-3dplot}

\usepgfplotslibrary{patchplots}

\usepackage{tikz}
\usepackage{tikzscale}

\usetikzlibrary{arrows,3d}

\usetikzlibrary{shapes.geometric}
\usetikzlibrary{positioning}
\usetikzlibrary{calc}
\usetikzlibrary{decorations.pathreplacing} 
\usetikzlibrary{decorations.markings}

\usetikzlibrary{spy}
\usetikzlibrary{snakes,arrows,shapes}    
\usetikzlibrary{spy}
\usetikzlibrary{backgrounds,tikzmark}  
\usetikzlibrary{decorations}
\usetikzlibrary{positioning}
\usetikzlibrary{patterns}
\usetikzlibrary{calc} 

\usepackage{booktabs}   
\usepackage{makecell} 
\usepackage{colortbl}   

\usepackage{xspace}
\usepackage{color}     
\usepackage{setspace}   
 
\usepackage{soul}
\usepackage{ulem}
\usepackage{currfile}
\usepackage{import}

\usepackage{authoraftertitle} 
\usepackage{authblk} 

\usepackage[noabbrev]{cleveref}

\usepackage{enumitem}
\usepackage{multicol} 
\usepackage{algorithm}
\usepackage{algpseudocode}
\usepackage{algorithmicx}
\usepackage[section]{placeins}

\usepackage{environ}

    \usepackage{framed}
    
    \setlist{noitemsep,topsep=2pt,parsep=0pt,partopsep=0pt}
  
    \usepackage{lmodern}

    \usepgfplotslibrary{external} 
    \usetikzlibrary{external}
   \pgfdeclarelayer{background}
   \pgfdeclarelayer{foreground}
   \pgfsetlayers{background,main,foreground}

\pgfplotsset{compat=1.8}

\newcommand{\findmax}[3]{
    \pgfplotstablesort[sort key={#2},sort cmp={float >}]{\sorted}{#1}%
    \pgfplotstablegetelem{0}{#2}\of{\sorted}%
    \let #3=\pgfplotsretval%
}

\pgfplotscreateplotcyclelist{myCycleListBlackAndWhite}
{%
  {black, mark=o},
  {black, mark=square},
  {black, mark=triangle},
  {black, mark=diamond}, 
  {black, mark=pentagon}, 
  {black, mark=*},
  {black, mark=square*},
  {black, mark=triangle*},
  {black, mark=diamond*}, 
  {black, mark=pentagon*}, 
}

\definecolor{darkgreen}{rgb}{0,0.4,0} 
\definecolor{darkbrown}{rgb}{0.5, 0.396, 0.09}

\definecolor{c1}{rgb}{0.0, 0.4196078431372549, 0.6431372549019608}
\definecolor{c2}{rgb}{1.0, 0.5019607843137255, 0.054901960784313725}
\definecolor{c3}{rgb}{0.6705882352941176, 0.6705882352941176,
0.6705882352941176} \definecolor{c}{rgb}{0.34901960784313724, 0.34901960784313724, 0.34901960784313724}
\definecolor{c4}{rgb}{0.37254901960784315, 0.6196078431372549,
0.8196078431372549} \definecolor{c}{rgb}{0.7843137254901961, 0.3215686274509804, 0.0}
\definecolor{c5}{rgb}{0.5372549019607843, 0.5372549019607843,
0.5372549019607843} \definecolor{c}{rgb}{0.6352941176470588, 0.7843137254901961, 0.9254901960784314}
\definecolor{c6}{rgb}{1.0, 0.7372549019607844, 0.4745098039215686}
\definecolor{c7}{rgb}{0.8117647058823529, 0.8117647058823529,
0.8117647058823529}

\pgfplotscreateplotcyclelist{myCycleListColor}
{%
  {blue, mark=o},
  {red, mark=square},
  {darkbrown, mark=triangle},
  {darkgreen, mark=diamond},
  {black, mark=pentagon},
  {blue, mark=*},
  {red, mark=square*},
  {darkbrown, mark=triangle*},
  {darkgreen, mark=diamond*},
  {black, mark=pentagon*},
}

\pgfplotsset{every axis/.append style= 
              {
                font=\scriptsize,
                mark size=3,
                legend style={font=\tiny, mark size=3, draw=none, fill=none},
                legend cell align=left,
                cycle list name=myCycleListColor,
                scaled y ticks = false,
				scaled x ticks = false,
				trim axis left,
				trim axis right,
				sharp plot,
				tick label style ={font=\tiny},
				label style ={font=\scriptsize},
				very thin,
				ymajorgrids=true,
				grid style=dotted,
				legend pos= north east,
				height=\figHeight, width=\figWidth,
              }
            }

\pgfplotsset{every tick label/.append style={font=\tiny}}
\pgfplotsset{every label/.append style={font=\scriptsize}}
\pgfplotsset{every axis legend/.append style={font=\tiny}}
\pgfplotsset{every axis plot/.append style={thick}}
            
\pgfplotstableset
{
  every odd row/.style={before row={\rowcolor[gray]{0.9}}},
  every head row/.style=
  {
    before row=
    {%
      \toprule
    },
    after row=\midrule
  },
  every last row/.style=
  {
    after row=\bottomrule
  }
}

\newif\ifdrawboundingbox

\pgfkeys
{ 
  /pgf/number format/.cd,
  fixed,
  set thousands separator={\,},
  1000 sep in fractionals,
}

\newcolumntype{x}[1]{>{\centering\arraybackslash}p{#1}}

\newcolumntype{C}[1]{>{\centering\arraybackslash}m{#1}}
\newcolumntype{R}[1]{>{\raggedright\arraybackslash}m{#1}}
\newcolumntype{L}[1]{>{\raggedleft\arraybackslash}m{#1}}

\newcommand\restr[2]{{
		\left.\kern-\nulldelimiterspace 
		#1 
		\right|_{#2} 
}}

\newcommand\compactDots{\ifmmode\ldots\else\makebox[10cm][c]{.\hfil.\hfil.}\fi}

\def\registered{{\ooalign{\hfil\raise .00ex\hbox{\tiny R}\hfil\cr\mathhexbox20D}}}

\makeatletter
\def\namedlabel#1#2{\begingroup
	\def\@currentlabel{#2}%
	\label{#1}\endgroup
}
\makeatother

\makeatletter
\renewcommand{\todo}[2][]{\tikzexternaldisable\@todo[#1]{#2}\tikzexternalenable}
\makeatother

\usepackage{etoolbox}
\AtBeginEnvironment{mysmallmatrix}{%
\setstretch{0.9}%
\setlength{\arraycolsep}{2pt}%
}%

{%
	\begin{matrix}%
}%
	{\end{matrix}}%

{%
	\begin{matrix}%
	}%
	{\end{matrix}}%
\AtBeginEnvironment{myverysmallmatrix}{%
	\footnotesize%
	\setstretch{0.55}%
	\setlength{\arraycolsep}{1.5pt}%
}%

\newcommand{\overbar}[1]{\mkern 0.1mu\overline{\mkern-0.1mu#1\mkern-0.1mu}\mkern 0.1mu}
\makeatletter
\def\bbordermatrix#1{\begingroup \m@th
	\@tempdima 4.75\p@
	\setbox\z@\vbox{%
		\def\cr{\crcr\noalign{\kern2\p@\global\let\cr\endline}}%
		\ialign{$##$\hfil\kern2\p@\kern\@tempdima&\thinspace\hfil$##$\hfil
			&&\quad\hfil$##$\hfil\crcr
			\omit\strut\hfil\crcr\noalign{\kern-\baselineskip}%
			#1\crcr\omit\strut\cr}}%
	\setbox\tw@\vbox{\unvcopy\z@\global\setbox\@ne\lastbox}%
	\setbox\tw@\hbox{\unhbox\@ne\unskip\global\setbox\@ne\lastbox}%
	\setbox\tw@\hbox{$\kern\wd\@ne\kern-\@tempdima\left[\kern-\wd\@ne
		\global\setbox\@ne\vbox{\box\@ne\kern2\p@}%
		\vcenter{\kern-\ht\@ne\unvbox\z@\kern-\baselineskip}\,\right]$}%
	\null\;\vbox{\kern\ht\@ne\box\tw@}\endgroup}
\makeatother

\title{Coupling of non-conforming trimmed isogeometric Kirchhoff-Love shells via a projected super-penalty approach.}

\author[1]{Luca Coradello\thanks{luca.coradello@epfl.ch, Corresponding Author}}
\author[3]{Josef Kiendl}
\author[1,2]{Annalisa Buffa}

\affil[1]{Chair of Numerical Modelling and Simulation,
 		  \'Ecole Polytechnique F\'ed\'erale de Lausanne, Lausanne, Switzerland.}

\affil[2]{Istituto di Matematica Applicata e Tecnologie Informatiche `E. Magenes' (CNR), Pavia, Italy.}
\affil[3]{Institute of Engineering Mechanics \& Structural Analysis,
	Bundeswehr University Munich, Germany}

\newcommand{\publicationDate}{\today}

\usepackage{fancyhdr}
\date{}
\fancypagestyle{plain}{%
  \fancyhf{}%
  \fancyfoot[L]{
  \vspace{-1.25cm} 
  \footnotesize{Preprint submitted to \textit{}  \hfill
  \publicationDate
  \\
  }} }
 
\begin{document}  
\normalem
\maketitle  
  
\vspace{-1.5cm} 
\hrule 
\section*{Abstract}

Penalty methods have proven to be particularly effective for achieving the required $C^1$-continuity in the context of multi-patch isogeometric Kirchhoff-Love shells. Due to their conceptual simplicity, these algorithms are readily applicable  to the displacement and rotational coupling  of trimmed, non-conforming surfaces. However, the accuracy of the resulting solution depends heavily on the choice of penalty parameters. Furthermore, the selection of these coefficients is generally problem-dependent and is based on a heuristic approach. Moreover, developing a penalty-like procedure that avoids interface locking while retaining optimal accuracy is still an open question.
This work focuses on these challenges. In particular, we devise a penalty-like strategy based on the $L^2$-projection of displacement and rotational coupling terms onto a degree-reduced spline space defined on the corresponding interface. Additionally, the penalty factors are completely defined by the problem setup and are constructed to ensure optimality of the method.
To demonstrate this, we asses the performance of the proposed numerical framework on a series of non-trimmed and trimmed multi-patch benchmarks discretized by non-conforming meshes. We systematically observe a significant gain of accuracy per degree-of-freedom and no interface locking phenomena compared to other penalty-like approaches. Lastly, we perform a static shell analysis of a complex engineering structure, namely the blade of a wind turbine.
\vspace{0.2cm} 
\hrule

\def\Estconst{3}
 
\vspace{0.25cm}
\noindent \textit{Keywords:} isogeometric analysis, multi-patch coupling, trimming, penalty method, Kirchhoff-Love shells. 
\vspace{0.25cm}
\hrule 


\section{Introduction}\label{sec:introduction}

The finite element method (FEM) is a well-established technology for the numerical simulation of a wide variety of engineering applications. It is also well-known that a significant amount of resources is invested into mesh generation and geometry clean-up~\citep{Cottrell2009}. To alleviate these burdens, isogeometric analysis (IGA) was introduced in the seminal paper~\citep{Hughes2005}. The main idea of IGA is to improve the interoperability between numerical simulations and Computer Aided Design (CAD) by employing the same mathematical objects used in the geometry description, namely B-splines and non-uniform rational B-splines (NURBS)~\citep{Piegl1995}, for the discretization of partial differential equations (PDEs). This shift in paradigm has paved the way for an extensive amount of research, where the reader is referred to~\citep{Hughes2005,Cottrell2009,Hughes20171} for a detailed review of the method and its recent state-of-the-art, while its mathematical analysis can be found~\citep{Bazilevs2006,Buffa2014}.

\noindent In particular, IGA has created a major impact on shells research. Historically, thick shell formulations of Reissner-Mindlin type are preferred in the finite element community~\citep{Bischoff2004} as they demand only $C^0$-continuity between elements. Classical formulations like the Kirchoff-Love shell are governed by fourth-order PDEs, resulting in a global $C^1$-continuity requirement which poses severe challenges to traditional finite element technologies. These obstacles are easily overcome within one isogeometric patch thanks to the higher continuity of B-splines, allowing to discretize higher-order PDEs directly in their primal form. We refer to~\citep{Kiendl2009,Kiendl2015,Kiendl2016} for a review of the method and several extensions in the scope of isogeometric Kirchhoff-Love shells, whereas other applications to high-order PDEs can be found in~\citep{Reali2015,Niiranen2017} for Kirchhoff plates, in~\citep{Gomez2008} for the Cahn-Hilliard equation, and in~\citep{BARTEZZAGHI2015} for the Laplace-Beltrami equation, respectively. Moreover, several other spline technologies have been successfully applied to the analysis of Kirchhoff-Love shells, for instance T-splines~\citep{BAZILEVS2012,Casquero2017,CASQUERO2020}, subdivision surfaces~\citep{Cirak2000}, PHT- and RHT-splines~\citep{NGUYENTHANH2011,NGUYENTHANH2017}, LR-splines~\citep{Proserpio2020}, and recently extended B-splines~\citep{schoellhammer2020consistent}.

However, in order to tackle structures of industrial relevance, two main issues need to be addressed. On one hand, the proper treatment of trimmed surfaces needs careful consideration, see~\citep{Marussig2018} for a review of the state-of-the-art and open challenges related to trimming. In this manuscript, we build a suitable high-order re-parametrization of those elements that are cut for integration purposes, by leveraging the tool presented in~\citep{Antolin2019}. This approach shares some similarities with the methods presented in~\citep{Breitenberger2015,Guo2018,Coradello2020Munich} for tackling trimmed shells and in~\citep{Kudela2015} in the scope of immersed methods.

On the other hand, complex geometries are typically described by multiple, non-conforming patches which, in turn, calls for a suitable coupling strategy to achieve the required $C^1$-continuity. Similarly to the nomenclature introduced in~\citep{Herrema2019}, we distinguish between $C^0$- or displacement continuity and $C^1$- or rotational continuity. The latter is not restricted to smooth interfaces but will also refer in the following to patches meeting at an arbitrary angle, which is preserved during deformation. This scenario is commonly encountered in complex engineering applications. 
In the literature, three methods are predominantly employed to achieve displacement and rotational continuity in a weak sense and they are briefly outlined in the following.

\noindent Mortar-type methods have been presented for patch coupling in~\citep{Horger2019,Hirschler2019} in the scope of Kirchhoff plates and Kirchhoff-Love shells, respectively, and have been generalized to arbitrary smoothness in~\citep{Dittmann2019}. It is well-known that mortar methods introduce additional artificial unknowns into the underlying system of equations to enforce the corresponding constraints, where the choice of discretization space for these Lagrange multipliers plays a pivotal role for the robustness of the method. In particular, the inf-sup stability is a crucial property, see~\citep{Boffi2013} for further details.

\noindent Another widespread coupling technique relies on Nitsche method, where the reader is referred to~\citep{Guo2015,Guo2018,NGUYENTHANH2017}. This family of methods is variationally-consistent and, generally speaking, it is more robust with respect to the choice of parameters compared to classical penalty approaches. However, Nitsche type algorithms are computationally less favorable as their implementation is problem-dependent and requires the computation of higher-order derivatives. In the case of shells, we highlight that these derivatives are defined on a manifold, which significantly increases the complexity.

\noindent Lastly, penalty-like methods are widely spread throughout numerous areas of engineering due to their versatility and ease of implementation. A variant thereof named bending strip was firstly studied in~\citep{Kiendl2010} in the scope of Kirchhoff-Love shells coupled along matching interfaces. Several other penalty-like approaches able to treat non-conforming discretization have followed, where the reader is referred to~\citep{Breitenberger2015,Duong2017}. However, the performance of these methods heavily relies on the choice of penalty parameters. It is well-known that these factors are problem-dependent and that finding suitable values is a tedious task, mainly based on a labor-intensive trial-and-error strategy. Moreover, finding a good balance between the imposition of the constraint, the condition number of the resulting system matrix and interface locking phenomena is of paramount importance for the accuracy of the resulting solution. Recently, these issues have attracted the attention of the research community. In~\citep{Herrema2019}, the parameters are scaled by the material properties, the thickness of the shell, the mesh size and by a single, user-defined, problem-independent factor which has been validated on an extensive series of benchmarks. In a similar manner,~\citep{Pasch2021} introduces an additional dependency of the penalty coefficients on various loading and boundary conditions. Concerning the interface locking, a possible remedy based on reduced integration has been proposed in~\citep{Leonetti2020}. Nevertheless, the development of a fully parameter-free penalty method which is insensitive to locking phenomena and retains optimal convergence is still in its infancy. 

Our contribution seeks to mitigate the aforementioned issues. Inspired by~\citep{Herrema2019} and taking~\citep{Coradello2020projected} as our starting point, we present an algorithm for achieving displacement and rotational continuity in trimmed shells which is inherently locking-free and where the penalty parameters are automatically defined by the problem setup, namely material properties, geometry of the structure and underlying discretization. Moreover, the aforementioned penalty factors are suitably built to retain the higher-order accuracy of B-splines. The methodology relies on the $L^2$-projection of the coupling terms along the associated interface onto a suitable degree-reduced space, where the stable $p/p-2$ pairing is employed~\citep{Brivadis2015}. We recall that the whole procedure is motivated by the analysis of the underlying perturbed saddle point problems, which gives us insight into the selection of appropriate parameters and into a proper way to eliminate the Lagrange multipliers associated to the constraints. Finally, we highlight that for splines of degree $p=2,3$ the projection is a local operation and therefore it introduces only a small overhead in the total run-time, making it computationally appealing.

\noindent Then, we verify numerically the robustness of the proposed coupling technique on various non-trimmed and trimmed examples and we compare it with other penalty-like methods. In all cases we observe an optimal convergence behavior, where no boundary locking effects are present. This yields a superior accuracy per degree-of-freedom (dof). Finally, we assess the applicability of our numerical framework to complex engineering structures. To demonstrate this, we perform a static analysis and simplified topology optimization of the DTU 10 MW Reference wind turbine blade~\citep{Bak2013}.

The paper is structured as follows. \Cref{sec:bsplines} provides to the reader the basic notation related to B-splines whereas \Cref{sec:trimming} introduces the fundamentals of trimming. \Cref{sec:formulation} explains in details the proposed method, focusing on how to choose the penalty parameters and on how to mitigate locking. In \Cref{sec:numericalExamples} the method is validated on a selection of non-trimmed and trimmed numerical benchmarks. Then, our approach is applied to the shell analysis of the DTU 10 MW Reference wind turbine blade, where a simplified topology optimization is performed on the stiffening webs. Finally, some conclusions are drawn in \Cref{sec:conclusions}.

\newcommand{\picsDir}{pictures/numericalExamples/pics}
\newcommand{\graphDir}{pictures/numericalExamples/graphs}
\newcommand{\dataDir}{pictures/numericalExamples/data}

\newtheorem{remark}{Remark}
\setlength{\extrarowheight}{0.05cm}

\section{B-splines in a nutshell} \label{sec:bsplines}
In this section, some definitions and fundamentals related to B-splines are reviewed. We refer the reader to~\citep{Piegl1995, Cottrell2009, Hoellig_book}, and references therein, for a comprehensive review of B-splines and NURBS and their role in isogeometric analysis.

\noindent Starting from two integers $p, n$, a univariate B-spline basis function $b_{i,p}$ of degree $p$ is generated starting from a non-decreasing sequence of real values referred to as knot vector, denoted in the following as $\Xi = \left\lbrace \xi_1, \ldots, \xi_{n+p+1} \right\rbrace$.
It is worth mentioning that the smoothness of the obtained B-spline basis is $C^{p-k}$ at every knot, where $k$ denotes the multiplicity of the considered knot, while it is $C^\infty$ elsewhere. In the remainder of this work, we consider only splines of maximum continuity, i.e. $C^{p-1}$, and degree $p \geq 2$.
The definition of multivariate B-splines $\mathcal{B}_{\mathbf{i},\mathbf{p}} (\bm{\eta} )$ is achieved in a straight-forward manner using the tensor product of univariate B-splines as:
\begin{align}
\mathcal{B}_{\mathbf{i},\mathbf{p}} (\bm{\eta} ) = \prod_{j=1}^{\widehat{d}} b_{i_{j},p_{j}}^{j} (\eta_j)\, ,
\end{align}
where $\widehat{d}$ denotes the dimension of the parameter space. Additionally, the multi-index $\mathbf{i}=\left\{ i_{1},...,i_{\widehat{d}}\right\} $
denotes the position in the tensor product structure and $\mathbf{p}=\left\{ p_{1},...,p_{\widehat{d}}\right\} $ indicates the vector of polynomial degrees, associated to the corresponding parametric dimension $\bm{\eta} = \eta_1,\ldots,\eta_{\widehat{d}}\,$, respectively. 
Finally, we denote by $\mathcal{Q}_0$ the Bezier mesh associated to the basis $\mathcal{B}_{\mathbf{i},\mathbf{p}}$.
Although omitted here, it is straightforward to extend the notation to NURBS, for details see~\citep{Cottrell2009}.
In the rest of the paper, without loss of generality, the degree vector $ \mathbf{p} $ will be considered equal in each parametric direction and therefore simplified to a single scalar value $p$. 


\section{Mathematical framework of trimming} \label{sec:trimming}
\newcommand{\bg}{}

\begin{figure}[!ht]
	\centering
	\includegraphics[width=0.6\textwidth]{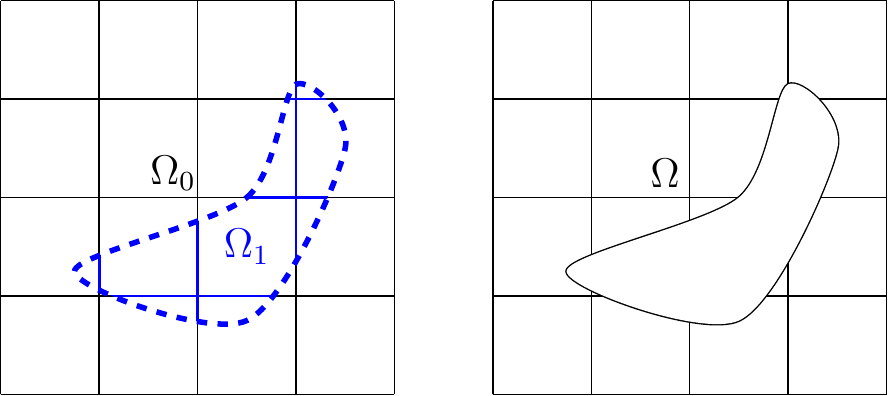}
	\caption{Example of a trimmed domain. From the untrimmed rectangular domain $\Omega_0$ the blue domain $\Omega_1$ is trimmed away, resulting in the final domain $\Omega$. For a correct interpretation of the colors, the reader is referred to the online version of this manuscript.} 
	\label{fig:trimming_exaple}
\end{figure}

In the following section, we summarize the basic mathematical foundation of isogeometric methods defined on trimmed domains, following closely the notation used in \citep{Antolin2019,Coradello2020}. For a detailed review of trimming and the current state-of-the-art in IGA we refer to \citep{Marussig2018} and references therein.
Let us define the domain $\Omega_0 \subset \mathbb{R}^d$, where $d$ is the dimension of the physical space of the problem at hand, described by a spline map $\mathbf{F}: \widehat{\Omega}_0 = [0, 1]^{\widehat{d}} \rightarrow \Omega_0$, where $\widehat{\Omega}_0$ denotes the parametric domain and we recall that $\widehat{d}$ represents its dimension. 
In particular, $\Omega_0$ is characterize as a linear combination of multivariate B-spline basis functions and corresponding control points as follows:
\begin{align}
	\Omega_0 = \mathbf{F}(\widehat{\Omega}_0) \quad \text{with} \quad \mathbf{F} (\bm{\eta} ) =\sum_{\mathbf{i}}\mathcal{B}_{\mathbf{i},p}(\bm{\eta} ) \mathbf{P}_{\mathbf{i}} \, .
\end{align}
Then, let us introduce $\Omega_1, \ldots, \Omega_N \subset \mathbb{R}^d$ Lipschitz-regular domains that define the trimming regions to be subtracted from $\Omega_0$. Consequently, the physical domain reads:
\begin{align} \label{eq:domain}
	\Omega = \Omega_0 \setminus \bigcup_{i = 1}^N \overline{\Omega}_i \, ,
\end{align}  
where an example is provided in~\Cref{fig:trimming_exaple} for the case $N=1$. 
We remark that the trimming operation does not change the underlying mathematical description of the original domain. Therefore, elements and associated basis functions are defined with respect to the untrimmed domain $\Omega_0$. 
\noindent Let us now introduce the B-spline basis of degree $p$ restricted to the corresponding parametric trimmed domain $\widehat{\Omega}$ as follows:
\begin{align}
\mathcal{B}_{\widehat{\Omega}} = \left\lbrace b \vert_{\widehat{\Omega}} : b \in \mathcal{B}_{\mathbf{i},p} \, \wedge \, \widehat{\Omega} \cap \text{supp}(b) \neq \varnothing \right\rbrace \, .
\end{align}
Similarly, we define the parametric mesh $\widehat{\mathcal{Q}}$ as the set of elements such that:
\begin{align}
	\widehat{\mathcal{Q}} = \left\lbrace Q \in \mathcal{Q}_0 : Q \cap \widehat{\Omega} \neq \varnothing \right\rbrace \, ,
\end{align}
where, in the following, we refer to $\mathcal{Q}$ as an active cell if $\mathcal{Q} \in \widehat{\mathcal{Q}}$. Consequently, the definition of physical mesh reads:
\begin{align}
	\mathcal{Q} = \left\lbrace \mathbf{F}(Q) \, : \,  Q \in \widehat{\mathcal{Q}} \right\rbrace \, ,
\end{align}
Finally, we introduce the approximation space formed by multivariate B-splines of degree $p$ restricted to a trimmed domain $\Omega$ as follows:
\begin{align}
S_{h}^p(\Omega) = \text{span }\left\lbrace b \circ \mathbf{F}^{-1} \, \vert \, b \in \mathcal{B}_{\widehat{\Omega}} \right\rbrace \, .
\end{align}

\section{The projected super-penalty method} \label{sec:formulation}

In this section, we extend the method studied by the authors in~\citep{Coradello2020projected} for coupling non-conforming Kirchhoff plates to the analysis of trimmed multi-patch Kirchhoff-Love shells. Motivated by the work presented in~\citep{Brivadis2015} in the context of isogeometric mortar methods, our strategy leverages the $L^2$ projection of the coupling terms at the interface, typically defined in terms of the degree $p$ of the solution space related to the corresponding patch, onto a reduced space of B-splines of degree $p - 2$ defined on the so-called active side of the interface. This procedure mitigates the locking phenomena due to the over-constraint of the solution space in the proximity of the corresponding coupling interface~\citep{Coradello2020Munich}. We remark that our method shares some similarities with the penalty coupling proposed in~\citep{Leonetti2020}.

\subsection{A review of differential geometry}

Let us review some fundamentals of differential geometry needed to describe the Kirchhoff-Love formulation, following closely the notation in~\citep{Benzaken2020}. 

\noindent Recalling the spline geometric mapping of the mid surface $\mathbf{F}: \widehat{\Omega} \rightarrow \Omega$,  the covariant basis vectors $\mathbf{a_\alpha}$ are defined as follows:
\begin{align}
\mathbf{a_\alpha}(\xi, \eta) = \mathbf{F}_{,\alpha}(\xi, \eta) \, ,
\end{align}
where the comma is used to indicate differentiation with respect to the corresponding curvilinear coordinate. Now, the unit normal vector to the mid surface of the shell $\mathbf{a}_3$ is computed as the normalized cross-product of the in-plane vectors $\mathbf{a}_\alpha$:
\begin{align}
\mathbf{a}_3 = \frac{\mathbf{a}_1 \times \mathbf{a}_2}{\vert \vert  \mathbf{a}_1 \times \mathbf{a}_2 \vert \vert} \, .
\end{align}
Then, let us introduce the covariant metric coefficients as:
\begin{align}
a_{\alpha \beta} = \mathbf{a_\alpha} \cdot \mathbf{a_\beta} \, .
\end{align}
Now, the contravariant basis vectors are defined via the following algebraic relationship:
\begin{align}
\mathbf{a^\alpha} \cdot \mathbf{a_\beta} = \delta^\alpha_\beta \, ,
\end{align}
where the corresponding covariant and contravariant metric coefficients are linked by the inverse operator:
\begin{align}
\Big[ a^{\alpha \beta} \Big] = \Big[ a_{\alpha \beta} \Big]^{-1} \, .
\end{align}
With these coefficients the contravariant basis vectors can be obtained as:
\begin{align}
\mathbf{a^\alpha} = a^{\alpha \gamma} \mathbf{a_\gamma} \, .
\end{align}
We now define the in-plane normal vector $\mathbf{n} = n^\alpha \mathbf{a}_\alpha$ to the boundary $\partial \Omega$, where $n^\alpha$ denotes its contravariant components. We highlight that $\mathbf{n}$ is contained in the tangent plane to the shell.
Finally, we define the transformation which maps Cartesian components to curvilinear ones as:
\begin{align}
\mathcal{Q}^i_\beta = \mathbf{e}^i \cdot \mathbf{a}_\beta \quad \text{and} \quad \mathcal{Q}^i_3 = \mathbf{e}^i \cdot \mathbf{a}_3 \, ,
\end{align}
where $\mathbf{e}^i$ represents the standard Euclidean basis.

\subsection{The weak form of the Kirchhoff-Love shell problem}
For the sake of conciseness, in the following we directly work in the discrete setting. The interested reader is referred to~\citep{Benzaken2020} for a rigorous derivation.
Let us consider as computational domain a manifold $\Omega \subset \mathbb{R}^3$ with a sufficiently smooth boundary $\Gamma = \partial \Omega$. Let us split the boundary $\Gamma = \partial \Omega$ into  a part associated to Dirichlet-type boundary conditions $\Gamma_D = \Gamma_\mathbf{u} \cup \Gamma_\mathbf{\theta}$ and a part corresponding to Neumann-type boundary conditions $\Gamma_N = \Gamma_\mathbf{T} \cup \Gamma_{B_{nn}}$ such that $\Gamma = \overline{\Gamma_D \cup \Gamma_N}$. It also holds that $\Gamma_\mathbf{u} \cap \Gamma_\mathbf{T} = \varnothing$ and $\Gamma_\mathbf{\theta} \cap \Gamma_{B_{nn}} = \varnothing$ due to the energetically conjugate nature of applied displacements and transverse shear, and applied rotations and bending moments, respectively.
Additionally, let us also introduce the set of corners as $\chi \subset \Gamma$, where this set can be further split into a Neumann part $\chi_N = \chi \cap \Gamma_N$ and a Dirichlet part $\chi_D = \chi \cap \Gamma_D$.
Let us also assume a given an applied body force $\tilde{\mathbf{f}} \in [ L^2(\Omega) ]^d$, a prescribed bending moment $ \tilde{B}_{nn} \in L^2(\Gamma_{B_{nn}})$, a prescribed Ersatz force $ \tilde{\mathbf{T}}\in L^2(\Gamma_{\mathbf{T}})$ as defined in~\citep{Benzaken2020} and a given twisting moment $\tilde{S} \in \mathbb{R}$ for all corners in $\chi_N$.
\noindent With these definitions at hand, the weak formulation of the Kirchhoff-Love shell reads, find $\mathbf{u}_h \in V_h$ such that:
\begin{align}\label{eq:weak_form_KL}
	a(\mathbf{u}_h, \mathbf{v}_h) = f(\mathbf{v}_h) \qquad \forall \mathbf{v}_h \in V_h\, ,
\end{align}
where the choice of discrete space $V_h \subset \left[S_{h}^p\right]^d$ depends in general on the boundary conditions of the problem at hand.
The bilinear form $a$ and linear form $f$ can be expanded, respectively, as follows:
\begin{align}
	a(\mathbf{u}_h, \mathbf{v}_h) &= \int_{\Omega} A(\mathbf{u}_h) \colon \alpha(\mathbf{v}_h) \textbf{d}\Omega +\int_{\Omega} B(\mathbf{u}_h) \colon \beta(\mathbf{v}_h) \textbf{d}\Omega  \nonumber \\
	f(\mathbf{v}_h) &= \int_{\Omega} \tilde{\mathbf{f}} \cdot  \mathbf{v}_h \textbf{d}\Omega + \int_{\Gamma_{\mathbf{T}}} \tilde{\mathbf{T}} \cdot  \mathbf{v}_h \textbf{d}\Gamma + \int_{\Gamma_{B_{nn}}} \tilde{B}_{nn}  \theta_n(\mathbf{v}_h) \textbf{d}\Gamma  + \sum_{e \in \chi_N} \left( \tilde{S} v_{3,h} \Big\vert_e \right)	\, ,
\end{align}
where we recall that the normal rotation $\theta_n(\mathbf{u})$ is given as:
\begin{align}
\theta_n(\mathbf{u}) = \mathcal{Q}^i_3 u_{i,\alpha} n^\alpha \, .
\end{align}
Then, by leveraging the in-plane projector $P = \mathbf{I} - \mathbf{a}_3 \otimes \mathbf{a}_3$, where $\mathbf{I}$ denotes the identity tensor, and the surface gradient $\nabla^{\star}$, the membrane and bending strain operators can be defined, respectively, as:
\begin{align}
	\alpha(\mathbf{v}_h) &= P \cdot \text{sym} \left( \nabla^{\star}(\mathbf{v}_h) \right) \cdot P \nonumber \\
	\beta(\mathbf{v}_h) &= -P \cdot \text{sym} \left( \mathbf{a}_3 \cdot \nabla^{\star} \nabla^{\star} (\mathbf{v}_h) \right) \cdot P \, ,
\end{align} 
where we highlight that the operator sym$(\cdot)$ returns the symmetric part of the input tensor. It is worth noting that the bending operator $\beta$ requires a global $C^1$-continuity of the basis to be well-defined. This is readily achieved within one patch by B-splines of degree $p \geq 2$.
Next, we can compute the corresponding stress operators by employing a constitutive law. In particular, if we consider a linear elastic model and we analytically integrate through the thickness, we can write:
\begin{align}
	A(\mathbf{u}_h) &= t \mathbb{C} \colon \alpha(\mathbf{u}_h) \nonumber \\
	B(\mathbf{u}_h) &= \frac{t^3}{12} \mathbb{C} \colon \beta(\mathbf{u}_h) \, ,
\end{align}
where the fourth-order tensor $\mathbb{C}$ for homogeneous materials can be expressed in curvilinear coordinates as:
\begin{align}
	\mathbb{C}=\mathbb{C}^{\alpha \beta \lambda \mu} \mathbf{a}_{\alpha} \otimes \mathbf{a}_{\beta} \otimes \mathbf{a}_{\lambda} \otimes \mathbf{a}_{\mu} \text { with } \mathbb{C}^{\alpha \beta \lambda \mu}=\frac{E}{2(1+v)}\left(a^{\alpha \lambda} a^{\beta \mu}+a^{\alpha \mu} a^{\beta \lambda}+\frac{2 v}{1-v} a^{\alpha \beta} a^{\lambda \mu}\right) \, ,
\end{align}
where $E$ and $\nu$ represent the Young's modulus and Poisson's ratio, respectively.

\begin{figure}[!ht]
	\centering
	\includegraphics[width=0.9\textwidth]{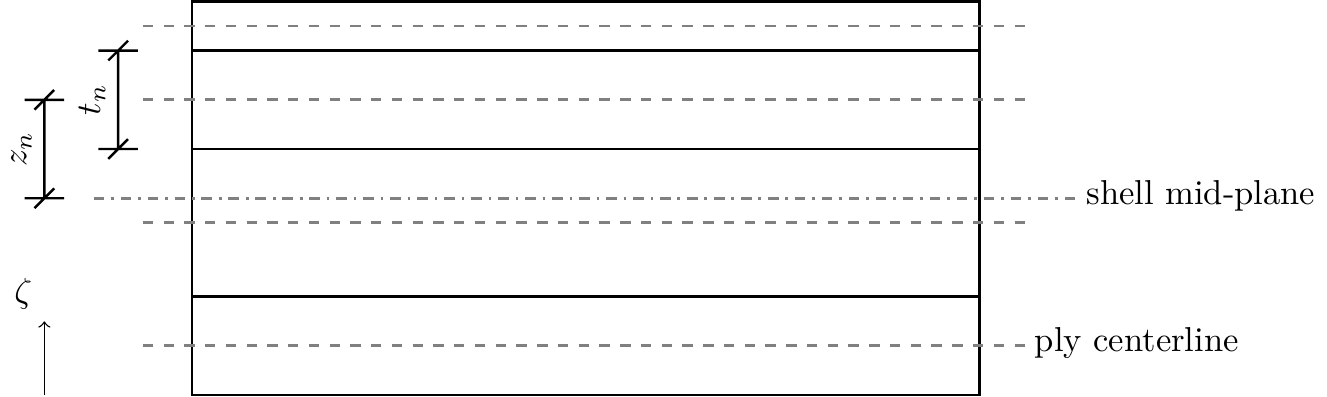}
	\caption{Example of a laminate along the thickness direction $\zeta$ formed by a non-uniform and non-symmetric ply sequence.} 
	\label{fig:composite_exaple}
\end{figure}

\noindent If we consider composite materials defined as a sequence of orthotropic plies, the bilinear form in~\eqref{eq:weak_form_KL} must be modified as explained in the following. 
Let us consider a stacking of plies, numbered by an index $n=1,\ldots,P$, where $P$ denotes the total number of plies. For each ply we can define the material tensor $\mathbb{C}_n$, obtained by transforming the corresponding orthotropic ply tensor from the local ply coordinates to the shell curvilinear reference frame, for further details see~\citep{KiendlThesis}.   
Now, following the classical theory of laminates~\citep{Reddy1999}, the homogenized extensional stiffness $\mathbb{A}$, the coupling stiffness $\mathbb{B}$ and the bending stiffness $\mathbb{D}$ are computed, respectively, as:
\begin{align}
	\mathbb{A} &= \int_{-t/2}^{t/2} \mathbb{C} \text{d}\zeta = \sum_{n=1}^{P} \mathbb{C}_n t_n \, , \nonumber \\
	\mathbb{B} &= \int_{-t/2}^{t/2} \zeta \, \mathbb{C} \text{d}\zeta = \sum_{n=1}^{P} \mathbb{C}_n t_n z_n \, , \nonumber \\
	\mathbb{D} &= \int_{-t/2}^{t/2} \zeta^2 \, \mathbb{C} \text{d}\zeta = \sum_{n=1}^{P} \mathbb{C}_n \left( t_n z_n^2 + \frac{t_n^3}{12} \right) \, ,
\end{align}
where $t_n$ indicates the thickness of the $n$-th ply and $z_n$ denotes the distance between the centroid of the $n$-th ply and the mid-plane of the shell, where an example is depicted in~\Cref{fig:composite_exaple}.

\noindent Then, the bilinear form associated to a laminate shell reads:
\begin{align}
	a(\mathbf{u}_h, \mathbf{v}_h) &= \int_{\Omega} \big( \mathbb{A} \colon \alpha(\mathbf{u}_h) +  \mathbb{B} \colon \beta(\mathbf{u}_h) \big) \colon \alpha(\mathbf{v}_h) \textbf{d}\Omega \nonumber \\
	&+ \int_{\Omega} \big( \mathbb{B} \colon \alpha(\mathbf{u}_h) +  \mathbb{D} \colon \beta(\mathbf{u}_h) \big) \colon \beta(\mathbf{v}_h) \textbf{d}\Omega \, ,
\end{align}
where for further details we refer to~\citep{KiendlThesis}.  
Finally,~\Cref{eq:weak_form_KL} can be summarized in matrix form as:
\begin{align}
	\boldsymbol{K} \overbar{\boldsymbol{u}} = \boldsymbol{f} \, ,
\end{align}
where $\boldsymbol{K}$ and $\boldsymbol{f}$ are denoted as the global stiffness matrix and force vector, respectively, and $\overbar{\boldsymbol{u}}$ represents the sought solution coefficients.
\subsection{The multi-patch setting}

\begin{figure}
	\centering
	\includegraphics[width=0.75\textwidth]{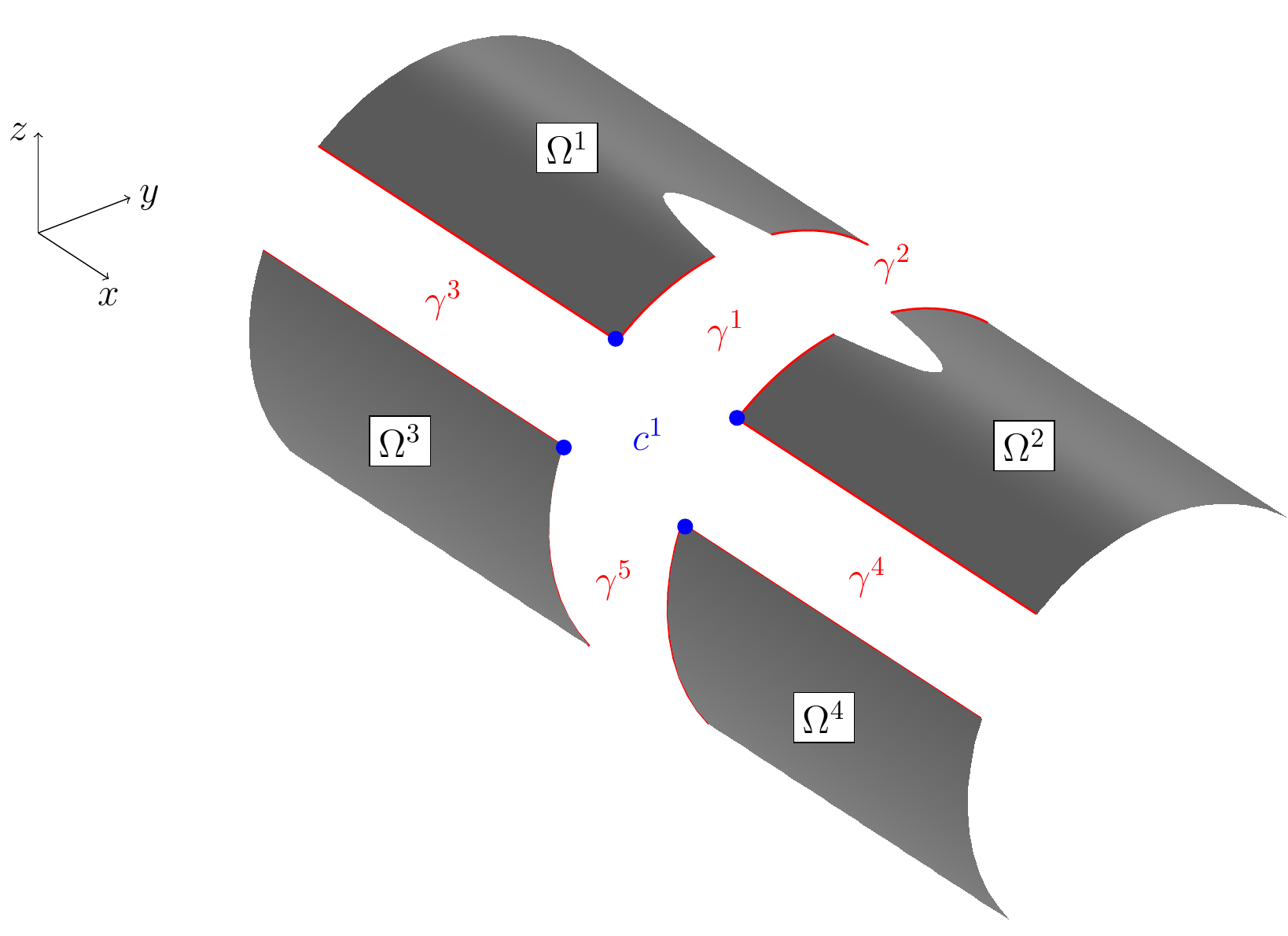}
	\caption{Example of four subdomains $\Omega^i, i=1,\ldots,4$ with their coupling interfaces $\gamma^{\ell}, \ell = 1,\ldots,5$, highlighted in red, and one corresponding cross-point $c^s, s=1$, represented by blue dots. Note that we have separated the subdomains for visualization purposes. For a correct interpretation of the colors, the reader is referred to the web version of this manuscript.}\label{fig:coupling_plate}
\end{figure}

Following closely the notation introduced in~\citep{Brivadis2015}, let us split the computational domain $\Omega$ into $N$ non-overlapping subdomains $\Omega^i$ such that:
\begin{align}
\overbar{\Omega} = \bigcup_{i=1}^N \overbar{\Omega^i} \, , \quad \text{where} \quad \Omega^i \cap \Omega^j = \varnothing \quad \text{for} \quad i \neq j \, .
\end{align}
In CAD terminology, $\Omega$ is a B-Rep, i.e. a collection of trimmed surfaces endowed with their topological information. In this work, similarly to~\citep{Breitenberger2015}, we use the so-called \textit{face-edge-vertex} B-Rep representation.
By leveraging the topology, we can then define the interface $\gamma^{\ell}$ between two adjacent trimmed patches $\Omega^{m}, \Omega^{n}, 1 \leq m, n \leq N$ as a common edge between their faces, see~\Cref{fig:coupling_plate} for an example on four patches. Note that two surfaces can share more than one edge.
Then, the skeleton $\Gamma$ is defined as the union of all common interfaces and reads:
\begin{align}
\Gamma = \bigcup_{\ell=1}^L \gamma^\ell \, ,
\end{align} 
where $L$ denotes the total number of interfaces and $\ell$ is an ordered index such that $1 \leq \ell \leq L$.
\begin{remark}
By a slight abuse of notation, $\gamma^\ell$ can represent both a trimmed or a non-trimmed coupling interface.
\end{remark}
\noindent Further, let us introduce the cross-points $c^{s}$ as the intersection of at least three shared edges and let us label them with an ordered index $c^{s}, s=1, \ldots, S$, see again~\Cref{fig:coupling_plate} for an illustration.
\begin{remark}
It is well-known that CAD softwares provide only an approximation of the true common edge $\gamma^\ell$ which depends on the chosen tolerance. For the sake of simplicity, in our derivation we assume exactness, or equivalently watertightness, of the geometric representation. 
From a computational standpoint, if the B-Rep is not watertight we perform a closest point projection of the relevant quantities, such as quadrature points and interface knots, onto the coupling edge. For further details we refer to~\citep{BAZILEVS2012}.
\end{remark}
\noindent Now, let us denote by $\mathbf{u}^m$ the value of the displacement field restricted to $\Omega^m$, and similarly $\mathbf{u}^n$ the value of the primary field on the neighboring subdomain $\Omega^n$. 
Then, for each interface $\gamma^{\ell}$ the following coupling conditions must be satisfied:
\begin{align}
\mathbf{u}^m - \mathbf{u}^n &= 0 \quad \text{on} \quad \gamma^{\ell}  \nonumber \\
\theta_n(\mathbf{u}^m) + \theta_n(\mathbf{u}^n) &= 0 \quad \text{on} \quad \gamma^{\ell} \, ,
\end{align}
which can be rewritten by leveraging the jump operator as:
\begin{align}\label{eq:coupling_conds}
\llbracket \mathbf{u} \rrbracket &= 0 \quad \text{on} \quad \gamma^{\ell} \nonumber \\
\llbracket \theta_n(\mathbf{u}) \rrbracket  &= 0 \quad \text{on} \quad \gamma^{\ell} \, . 
\end{align}
\subsection{The projected super-penalty formulation}
\begin{figure}
	\centering
	\begin{subfigure}[t]{0.7\textwidth}
		\centering
		\includegraphics[width=\textwidth]{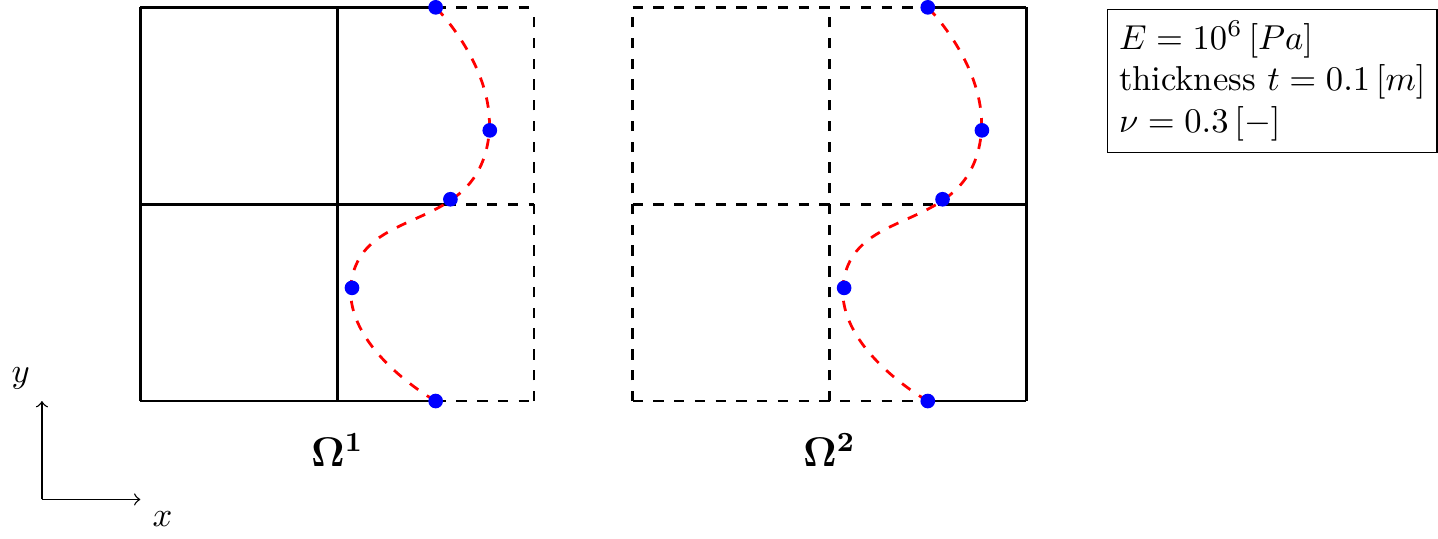}
		\caption{Initial discretization. The trimmed interface is represented by the dashed red line. Blue dots denote the location of the B-spline curve knots.}
	\end{subfigure}	
	\begin{subfigure}[t]{0.495\textwidth}
		\centering
		\input{\graphDir/convergence_shell_knotsOff_conforming.tex}
		\caption{Error $H^2$, without knots.}
	\end{subfigure}
	\hfill
	\begin{subfigure}[t]{0.495\textwidth}
		\centering
		\input{\graphDir/convergence_shell_knotsOn_conforming.tex}
		\caption{Error $H^2$, with knots.}
	\end{subfigure}
	\caption{Influence of internal knots of the coupling curve on a two trimmed plates example.}
	\label{fig:influence_knot}
\end{figure}

\begin{figure}
	\centering
	\includegraphics[width=0.65\textwidth]{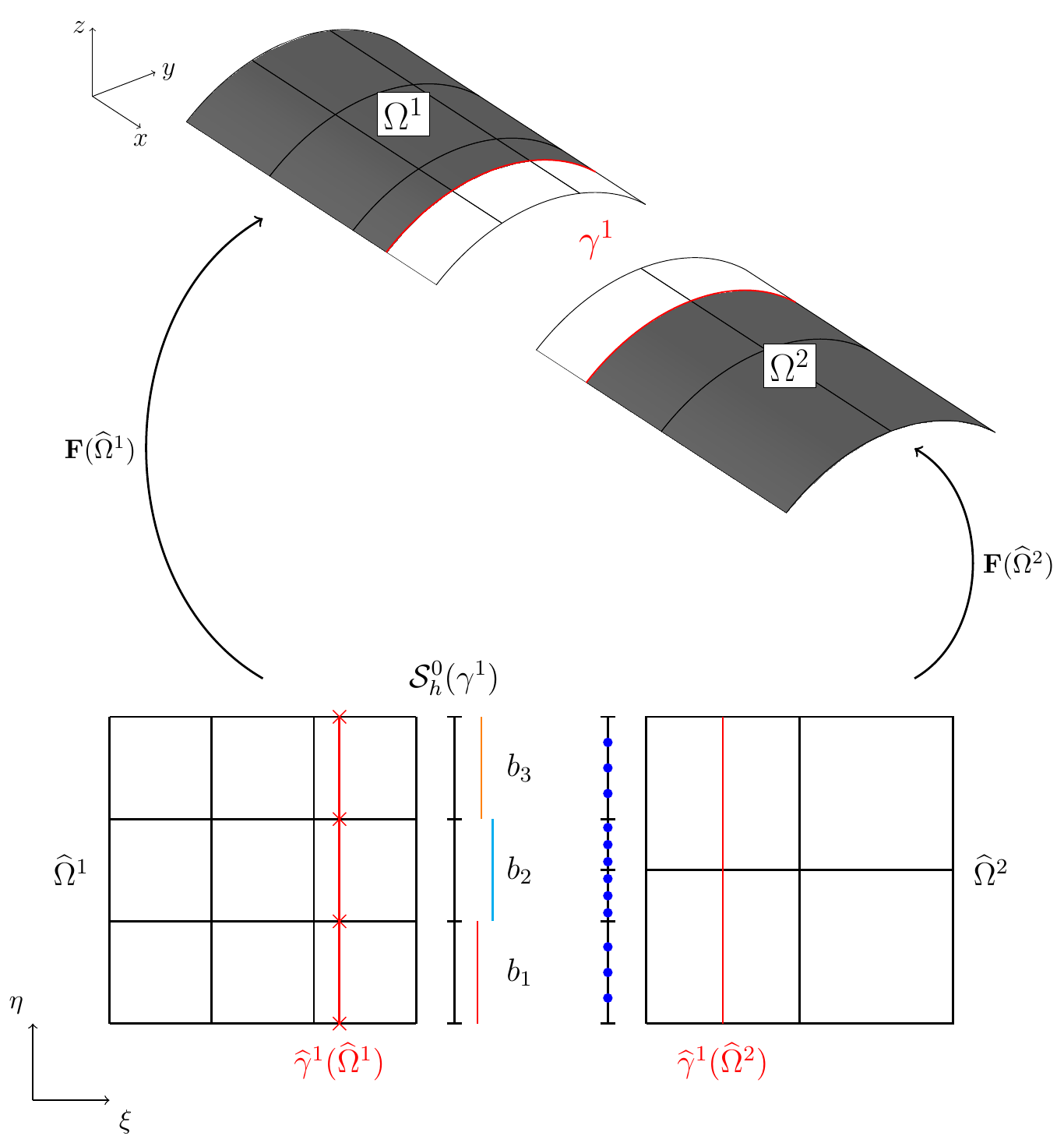}
	\caption{Example of the projection setup on a coupling interface for B-splines of degree $p = 2$. We arbitrarily select the finer mesh (on $\Omega^1$ in this example) to define the projection space ${\mathcal{S}^{p-2}_{h}(\gamma^{\ell})}$, where the intersections between the parametric coupling curve $\widehat{\gamma}^1(\widehat{\Omega}^1)$ and the knot lines of $\Omega^1$ are represented by red crosses. Additionally, an intersection mesh at the interface is created only for integration purposes to properly compute the projected penalty terms in~\Cref{eq:weak_form_with_penalty}. The $p+1$ integration points are schematically represented by blue dots. Note that we have separated the subdomains for visualization purposes.}\label{fig:proj_setup}
\end{figure}

Following the notation presented in~\citep{Coradello2020projected}, let us introduce for each patch $\Omega^i$ the following space:
\begin{align}
Z_{i,h}=\mathrm{span} \{ b \in \left[\mathcal{S}^p_{h}(\Omega^i) \right]^d \} \, .
\end{align} 
Additionally, we denote by $V_{i,h} \subset Z_{i,h}$ the finite-dimensional space given by the span of splines associated to subdomain $\Omega^i$, where the exact definition of $V_{i,h}$ depends on the set of boundary conditions of the problem at hand.
This allows us to introduce the following finite-dimensional space,
\begin{align}
V_h &= \left\lbrace v \in L^2(\Omega) \, \vert \, v \in V_{i,h} \: \forall i=1,\ldots,N \:\text{ and } v \text{ is continuous in } c^s, \:s=1,\ldots,S \right\rbrace \, ,
\end{align}
where we highlight the $C^0$-continuity requirement at the cross-points $c^s$.
Furthermore, for each interface $\gamma^{\ell}$, we introduce the associated knot vector $\Xi^{\ell}$. The latter is constructed as follows. First, we arbitrarily choose one of the neighboring patches as active. Then, we build $\Xi^{\ell}$ by intersecting the knot lines of the active patch and $\gamma^{\ell}$. We highlight that this operation can be performed directly in the parameter space of the active surface, since the B-Rep structure provides a representation of $\gamma^{\ell}$ in the parameter space of both surfaces. For each surface, we denote the latter representation by $\widehat{\gamma}^{\ell}(\widehat{\Omega}^i)$, see~\Cref{fig:proj_setup} for an example.
\begin{remark}
{\color{black}At this stage, in the spirit of developing a simple and efficient method, we disregard the internal knots of the coupling curve for the construction of $\Xi^{\ell}$. We highlight that the number of these knots depends on the chosen tolerance in the CAD model, with this number being potentially large. We are aware that this choice could potentially yield a loss of optimality of the method, but for smooth interfaces this effect is negligible. We verify this numerically on two trimmed patches in~\Cref{fig:influence_knot}, coupled along a $C^1$-continuous quadratic B-spline curve. Indeed, for $p=2,3$ the results are practically indistinguishable, whereas only minor differences are present for the case $p=4$. Although outside the scope of this work, finding a simple way to remove this source of sub-optimality constitutes a future research direction.}
\end{remark}	
Then, we build the isogeometric space $\mathcal{S}^{p-2}_{h}(\gamma^{\ell})$ leveraging the $p/p-2$ pairing. Assuming B-splines of maximum smoothness, this space is obtained by removing from $\Xi^{\ell}$ the first and last two knots, where an illustrative example is given in~\Cref{fig:proj_setup} for bivariate B-splines of degree $p=2$ and corresponding $p-2 = 0$ degree-reduced splines defined on the interface knot vector \mbox{$\Xi^{\ell} = \left[ 0 \,\, 1/3 \,\, 2/3 \,\,  1 \right]$}.
\begin{remark}
The $p/p-2$ pairing has been proven to be inf-sup stable in the context of isogeometric mortar methods in~\citep{Brivadis2015} and it has been extended to the coupling of non-trimmed Kirchhoff plates in~\citep{Coradello2020projected}. Although its stability for trimmed geometries has not been rigorously studied, we verify numerically its applicability to the coupling of trimmed Kirchhoff-Love shells.
\end{remark}
Consequently, let us define the following space:
\begin{align}
	Q_h =  \left\lbrace \mu \in L^2(\Gamma) \, \vert \, \mu \in \mathcal{S}^{p-2}_{h}(\gamma^{\ell}) \: \forall \ell=1,\ldots,L \right\rbrace \, , 
\end{align} 
which is used to characterize the Lagrange multipliers associated to the coupling conditions.
\noindent We are now ready to define the discretized version of the singularly-perturbed saddle point problem associated to the Kirchhoff-Love shell. Without loss of generality, let us consider homogeneous Neumann-type boundary conditions. Then, the saddle problem reads: find $\left( \mathbf{u}_h, \bm{\lambda}_{1,h}, \lambda_{2,h}  \right) \in V_h \times [Q_h]^d \times Q_h$ such that:
\begin{align}\label{eq:discrete_perturbed_problem}
\sum_{i=1}^{N} \left(\int_{\Omega^i} A(\mathbf{u}_h) \colon \alpha(\mathbf{v}_h) + B(\mathbf{u}_h) \colon \beta(\mathbf{v}_h) \right) + \sum_{\ell=1}^{L}\left(\int_{\gamma^{\ell}} \llbracket \mathbf{v}_h \rrbracket \bm{\lambda}_{1,h} + \int_{\gamma^{\ell}} \left\llbracket \theta_n(\mathbf{v}_h) \right\rrbracket \lambda_{2,h}\right) &= (\mathbf{f},\mathbf{v}_h) & \: &\forall \mathbf{v}_h \in V_h  \nonumber \\
\sum_{\ell=1}^{L} \left( \int_{\gamma^{\ell}} \llbracket \mathbf{u}_h \rrbracket \bm{\mu}_{1,h} - \frac{1}{\alpha_{\text{disp}}^\ell} \int_{\gamma^{\ell}} \bm{\lambda}_{1,h} \bm{\mu}_{1,h} \right) &= 0  & \quad &\forall \bm{\mu}_{1,h} \in [Q_h]^d \nonumber \\
\sum_{\ell=1}^{L} \left( \int_{\gamma^{\ell}} \left\llbracket \theta_n(\mathbf{u}_h) \right\rrbracket \mu_{2,h} - \frac{1}{\alpha_{\text{rot}}^\ell} \int_{\gamma^{\ell}} \lambda_{2,h} \mu_{2,h} \right) &=0 & \quad &\forall \mu_{2,h} \in Q_h \, , 
\end{align}
where we have introduced the parameters $\alpha_{\text{disp}}^\ell$ and $\alpha_{\text{rot}}^\ell$ corresponding to the displacements and normal rotations, respectively. For a rigorous derivation of the singularly-perturbed saddle point formulation in the scope of Kirchhoff plates we refer to~\citep{Coradello2020projected}.
As highlighted in~\citep{Herrema2019,Pasch2021}, these coefficients depend in general on the problem definition, e.g. the material parameters, the thickness of the shell, the applied boundary conditions, the mesh size and discretization degree, where a precise definition of our parameters will be provided in a later section.
Let us now eliminate the Lagrange multipliers and rewrite~\eqref{eq:discrete_perturbed_problem} only in terms of the displacement field. In particular, rearranging the second and third equations we obtain:
\begin{align}\label{eq:lagrange_multipliers}
{\bm{\lambda}_{1,h}}|_{\gamma^{\ell}} &= \alpha_{\text{disp}}^\ell \Pi^{\ell} \llbracket \mathbf{u}_h \rrbracket \nonumber \\ 
{\lambda_{2,h}}|_{\gamma^{\ell}} &= \alpha_{\text{rot}}^\ell \Pi^{\ell} \left\llbracket \theta_n(\mathbf{u}_h) \right\rrbracket \, , 
\end{align}
where, with a slight abuse of notation, $\Pi^{\ell}$ stands for the $L^2$-projection, defined on the interface $\gamma^{\ell}$, onto the degree-reduced space $\left[\mathcal{S}^{p-2}_{h}(\gamma^{\ell})\right]^d$ related to the displacements and onto the space $\mathcal{S}^{p-2}_{h}(\gamma^{\ell})$ associated to the normal rotations, respectively.

\noindent By substituting~\Cref{eq:lagrange_multipliers} into the first line of~\eqref{eq:discrete_perturbed_problem} and leveraging the properties of the $L^2$-projection, we obtain:
\begin{align}\label{eq:weak_form_with_penalty}
	a_{p}(\mathbf{u}_h, \mathbf{v}_h) &= \sum_{i=1}^{N} a^{(i)}(\mathbf{u}_h, \mathbf{v}_h) + \sum_{\ell=1}^{L}  \left( \int_{\gamma^{\ell}} \alpha_{\text{disp}}^{(\ell)} \Pi^{\ell} \llbracket \mathbf{v}_h \rrbracket \cdot \Pi^{\ell} \llbracket \mathbf{u}_h \rrbracket \, + \int_{\gamma^{\ell}} \alpha_{\text{rot}}^{(\ell)} \Pi^{\ell} \llbracket \theta_n(\mathbf{v}_h) \rrbracket \, \Pi^{\ell} \llbracket \theta_n(\mathbf{u}_h) \rrbracket \right) \, .
\end{align}
These coupling terms weakly impose the transmission conditions in~\eqref{eq:coupling_conds} on the displacements and normal rotations, respectively.
\begin{remark}
	From a computational standpoint, we rewrite the coupling term associated to the rotations in~\eqref{eq:weak_form_with_penalty} as defined in~\citep{Herrema2019}, where the constraint is recast into two complementary terms. This ensures a non-zero penalty contribution for patches meeting at an arbitrary angle. Then, the $L^2$-projection of these terms is performed. For further details, we refer to~\citep{Herrema2019} and references therein.
\end{remark}
\noindent Now, let us further characterize the aforementioned projection from a computational viewpoint. Let us consider a generic function $u \in V_h$ defined as the linear combination of basis functions and their corresponding coefficients $\boldsymbol{\hat{u}}$ as:
\begin{align}\label{eq:coefs1}
u = \sum_{i} \mathcal{B}_{i} \hat{u}_i \quad i = 1,\ldots,\text{dim} (V_h) \, .
\end{align}
Similarly, its projection $\Pi^{\ell}(u)$ onto the space $\mathcal{S}^{p-2}_{h}(\gamma^{\ell})$ can be written as another linear combination of spline functions and their associated coefficients $\boldsymbol{\tilde{u}}$:
\begin{align}\label{eq:coefs2}
\Pi^{\ell}(u) = \sum_{j} b_{j} \tilde{u}_j \quad j = 1,\ldots,\text{dim} (\mathcal{S}^{p-2}_{h}(\gamma^{\ell})) \, .
\end{align}
The orthogonality of the projection can now be expressed as:
\begin{align}\label{eq:opt_proj}
\int_{\gamma^{\ell}} \Pi^{\ell}(u) b = \int_{\gamma^{\ell}} u b \quad \forall b \in \mathcal{S}^{p-2}_{h}(\gamma^{\ell}) \, ,
\end{align}
which can be rewritten in matrix form by substituting~\Cref{eq:coefs1,eq:coefs2} into~\Cref{eq:opt_proj} as follows:
\begin{align}\label{eq:proj_sys}
\mathcal{M} \boldsymbol{\tilde{u}} = \mathcal{F} \boldsymbol{\hat{u}} \, ,
\end{align}
where $\mathcal{M}$ denotes the mass matrix associated to the degree-reduced basis and $\mathcal{F}$ represents the right-hand-side matrix corresponding to the inner product between the basis functions in $\mathcal{S}^{p-2}_{h}(\gamma^{\ell})$ and $V_h$, respectively. In particular, for the projection of the displacement term introduced in~\Cref{eq:lagrange_multipliers}, $\mathcal{F}_{\text{disp}}$ is defined as the inner product between the splines in $\left[\mathcal{S}^{p-2}_{h}(\gamma^{\ell})\right]^d$ and the jump of the basis functions in $V_h$. Analogously for the rotational term, $\mathcal{F}_{\text{rot}}$ is assembled as the inner product between the basis functions in $\mathcal{S}^{p-2}_{h}(\gamma^{\ell})$ and the jump of discrete normal rotations in $V_h$. Similarly, we distinguish between the mass matrix $\mathcal{M}$ associated to the splines in $\mathcal{S}^{p-2}_{h}(\gamma^{\ell})$ and its vectorial counterpart $\boldsymbol{\mathcal{M}}$ corresponding to the functions in $\left[\mathcal{S}^{p-2}_{h}(\gamma^{\ell})\right]^d$.
\noindent With these definitions at hand, we summarize the computation of the projected terms in~\Cref{alg:proj}.
\begin{algorithm}[H] 
	\begin{algorithmic}[1]
		\Procedure{Computation of the penalty terms}{} 
		\For{\textbf{each} interface $\gamma^{\ell}$ in $\Gamma$}
			\State Build the spaces $\mathcal{S}^{p-2}_{h}(\gamma^{\ell})$ and $\left[\mathcal{S}^{p-2}_{h}(\gamma^{\ell})\right]^d$
			\State Build the intersection mesh for integration
			\State $\boldsymbol{\tilde{u}}_{\text{disp}}\leftarrow$ solve~\Cref{eq:proj_sys} with $\boldsymbol{\mathcal{M}}$ and $\mathcal{F}_{\text{disp}}$
			\State $\boldsymbol{\tilde{u}}_{\text{rot}}\leftarrow$ solve~\Cref{eq:proj_sys} with $\mathcal{M}$ and $\mathcal{F}_{\text{rot}}$
			\State $\boldsymbol{K}$ = $\boldsymbol{K}$ + $\alpha_{\text{disp}}^{(\ell)} \boldsymbol{\tilde{u}}_{\text{disp}}^\top \boldsymbol{\mathcal{M}} \boldsymbol{\tilde{u}}_{\text{disp}}$
			\State $\boldsymbol{K}$ = $\boldsymbol{K}$ + $\alpha_{\text{rot}}^{(\ell)} \boldsymbol{\tilde{u}}_{\text{rot}}^\top \mathcal{M} \boldsymbol{\tilde{u}}_{\text{rot}}$
		\EndFor
		\EndProcedure
	\end{algorithmic} 
	\caption{Computation of the penalty terms in~\Cref{eq:weak_form_with_penalty}.}\label{alg:proj}
\end{algorithm}    
\noindent Lastly, we remark that the solution of~\Cref{eq:proj_sys} is computationally inexpensive for B-splines of degree $p=2,3$ associated to a reduced space of degree $p-2=0,1$, respectively, for which the mass matrix is either diagonal or can be lumped.

\subsubsection{Selection of penalty parameters}

It is well-known that the perturbed problem~\eqref{eq:discrete_perturbed_problem} is variationally consistent only if we select $\alpha_{\text{disp}}^\ell = \alpha_{\text{rot}}^\ell \rightarrow \infty \, \, \ell = 1,\ldots,L$. 
However, the well-posedness of the underlying problem is insensitive to the choice of the parameters $\alpha_{\text{disp}}^\ell$ and $\alpha_{\text{rot}}^\ell$. Therefore, our method is inherently free from boundary locking, independently of the choice of penalty values, see~\citep{Coradello2020projected} for further details in the context of isogeometric Kirchhoff plates. This allows us to select $\alpha_{\text{disp}}^\ell$ and $\alpha_{\text{rot}}^\ell$ to guarantee the high-order convergence rates achievable by B-splines.
Furthermore, in the spirit of developing a parameter-free penalty method, we modify the choice proposed in~\citep{Herrema2019}, scaling the displacement and rotation penalty parameters by the physical constants of the underlying problem, the local mesh size, the spline degree and the geometry. For homogeneous isotropic materials they read:
\begin{align}\label{eq:projectedTerms}
\alpha_{\text{disp}}^\ell &= (\vert \gamma^{\ell} \vert )^{\beta-1} \frac{E t}{(h_\ell)^\beta (1 - \nu^2)}  \nonumber\\
\alpha_{\text{rot}}^\ell &= (\vert \gamma^{\ell} \vert)^{\beta-1} \frac{E t^3}{12 (h_\ell)^\beta (1 - \nu^2)} \, ,
\end{align} 
where the measure of $\gamma^{\ell}$ serves as a characteristic length and the exponent $\beta$ is chosen solely to ensure the optimal convergence of the method. Therefore, it must be a function of the degree $p$ of the underlying discretization.
Numerically we have observed that the scaling factor $\beta = p-1$ in~\eqref{eq:projectedTerms} is necessary to attain optimal convergence of the method in the $H^2$ norm, whereas for a scaling of $\beta = p$ we noticed optimality in the $H^2$ and $H^1$ norms. Finally, a factor of $\beta = p+1$
provides optimality in the $H^2$, $H^1$ and $L^2$ norms. If not stated otherwise, we will use $\beta = p+1$ in all our numerical examples.
In case of orthotropic laminates, we adapt the minimum strategy presented in~\citep{Herrema2019}, where the minimum local stiffness between adjacent patches $\Omega^m$ and $\Omega^n$ is used. Consequently, the penalty parameters are defined as:
\begin{align}\label{eq:projectedTerms_composite}
	\alpha_{\text{disp}}^\ell &= (\vert \gamma^{\ell} \vert )^{\beta-1} \frac{\min({\max_{i,j}(\mathbb{A}_{ij}^{(m)}),\, \max_{i,j}(\mathbb{A}_{ij}^{(n)})})}{(h_\ell)^\beta}  \nonumber\\
	\alpha_{\text{rot}}^\ell &= (\vert \gamma^{\ell} \vert)^{\beta-1} \frac{\min({\max_{i,j}(\mathbb{D}_{ij}^{(m)}), \,\max_{i,j}(\mathbb{D}_{ij}^{(n)})})}{(h_\ell)^\beta} \, .
\end{align} 
Note that all of these parameters are known and depend only on the problem definition, meaning that no user-defined factor is required. Moreover, it is straightforward to check that the penalty terms are dimensionally consistent with respect to their corresponding energy contribution in the weak form~\eqref{eq:weak_form_with_penalty}. 
\begin{remark}
Clearly, the choice of $\beta$ influences the condition number of the associated system matrix. This, together with small trimmed elements, can potentially yield ill-conditioned systems of equations and, consequently, loss of accuracy due to numerical round-off errors. In the context of trimmed single-patch shells, a possible remedy based on extended B-splines has been studied in~\citep{schoellhammer2020consistent}. Furthermore, in the scope of immersed methods, an ad-hoc multigrid preconditioner has been developed in~\citep{de_prenter_multigrid_2020}.
In this contribution, we employ a direct solver where the stiffness matrix is preconditioned by a simple diagonal scaling. This seems to suffice for the level of accuracy reached in our numerical experiments. We remark that a thorough study of the condition number in the context of trimmed multi-patch Kirchhoff-Love shells is beyond the scope of this paper.
\end{remark}
\renewcommand{\graphDir}{pictures/numericalExamples/graphs}
\renewcommand{\dataDir}{pictures/numericalExamples/data}
\renewcommand{\picsDir}{pictures/numericalExamples/pics}

\section{Numerical Examples}
\label{sec:numericalExamples}

In this section we assess the performance of the proposed coupling technique with several numerical examples defined both on trimmed and untrimmed, non-conforming, multi-patch geometries. 
All the numerical examples presented in the following have been implemented in the open-source and free Octave/Matlab package \textit{GeoPDEs}~\citep{Vazquez2016}, where the reparametrization of the trimmed elements for integration purposes is provided by the tool presented in~\citep{Antolin2019}.	
The analytical shell solutions are taken from the new shell obstacle course studied in~\citep{Benzaken2020}, where the exact manufactured functions are evaluated in the freely-available Mathematica notebook\footnote{https://github.com/wdas/shell-obstacle-course} with 100 digits of precision. Moreover, similarly to~\citep{Benzaken2020}, for every element we employ $25\times25$ quadrature points to properly capture the highly non-linearity of the quantities of interest. The results of these computations are then imported into \textit{GeoPDEs}.
Also, in all examples taken from~\citep{Benzaken2020}, we derive from the manufactured solution and apply on the entire boundary $\partial \Omega$ non-homogeneous Dirichlet boundary conditions for the displacements and non-homogeneous Neumann boundary conditions for the bending moments.

\noindent Finally, throughout this section, we compare our choice of penalty factors to a classical approach where the parameters are kept constant:   
\begin{align}
	\alpha_{\text{disp}}^\ell &= 10^3 E \nonumber\\
	\alpha_{\text{rot}}^\ell &= 10^3 E \, ,
\end{align} 
and to the method proposed in~\citep{Herrema2019}:
\begin{align}\label{eq:projectedTerms}
	\alpha_{\text{disp}}^\ell &= \delta \frac{E t}{(h_\ell) (1 - \nu^2)}  \nonumber\\
	\alpha_{\text{rot}}^\ell &= \delta \frac{E t^3}{12 (h_\ell) (1 - \nu^2)} \, ,
\end{align} 
where the problem-independent, user-defined parameter $\delta = 10^3$ has been numerically validated on an extensive series of benchmarks.

\subsection{Four non-trimmed planar patches}
\begin{figure}
	\centering
	\begin{subfigure}[t]{0.495\textwidth}
		\centering
		\includegraphics[width=0.95\textwidth]{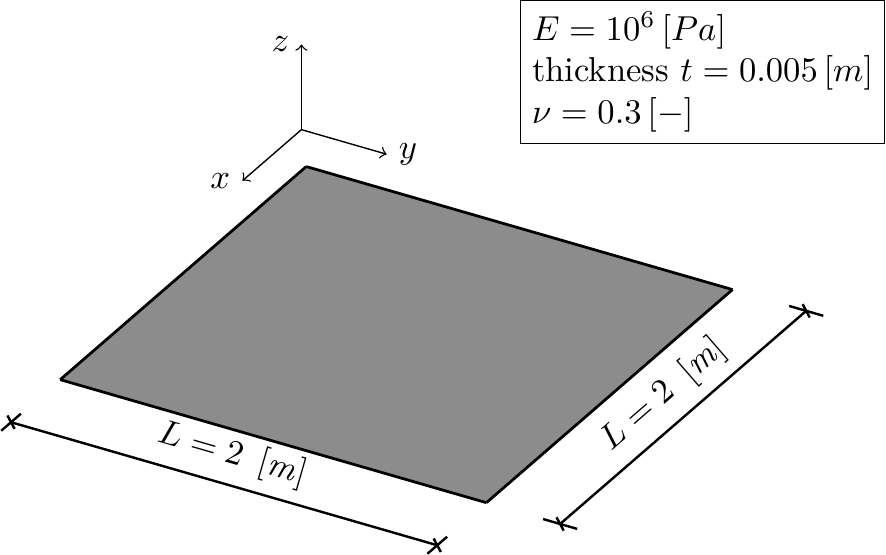}
		\caption{Geometry setup and physical parameters.}
	\end{subfigure}
	\hfill
	\begin{subfigure}[t]{0.495\textwidth}
		\centering
		\includegraphics[width=0.85\textwidth]{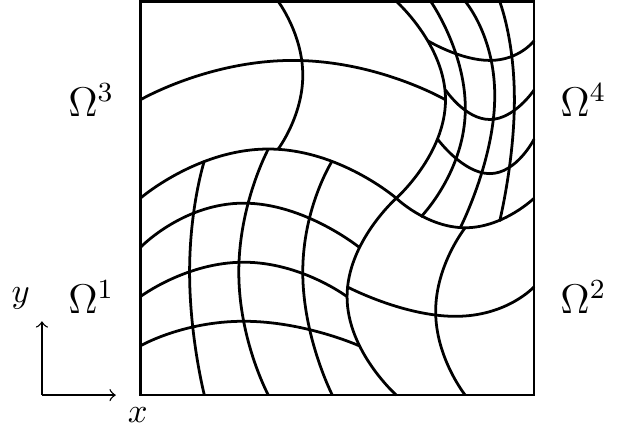}
		\caption{Initial discretization.}
	\end{subfigure}	
	\caption{Problem setup and initial non-trimmed, non-conforming, multi-patch discretization for the four planar patches example.}
	\label{fig:setup_curved_4patches}
\end{figure}

The first example is meant to test and verify the implementation of our strategy in a non-trimmed planar setting. The geometrical setup is taken from~\citep{Coradello2020projected}, where this problem is studied in the context of Kirchhoff plates. In particular, the domain $\Omega = [0,2] \times [0,2]$ is subdivided into four non-conforming patches $\Omega^i, i=1,\ldots,4$ coupled along curved interfaces, see~\Cref{fig:setup_curved_4patches}. To enforce the non-conformity of the latter, the initial interface knots have been shifted by the irrational factor $\sqrt{2}/100$. In our problem definition, we set the Young's modulus $E = 10^6 \, [Pa]$, the thickness of the plate $t=0.005 \, [m]$ and the Poisson's ratio $\nu = 0.3 \, [-]$, respectively.
Then, to verify the theoretical orders of convergence, we compute the approximation error in the $L^2$ and $H^2$ norms with respect to a manufactured smooth solution of the form:
\begin{align}
\bm{u}^{\text{ex}} (x, y, z) = 
\begin{pmatrix}
	u^x \\
	u^y \\
	u^z 
\end{pmatrix} = 
\begin{pmatrix}
	\sin(\pi  x) \sin(\pi  y) \\
	\sin(\pi  x) \sin(\pi  y) \\
	\sin(\pi  x) \sin(\pi  y) 
\end{pmatrix} 
\, .
\end{align} 

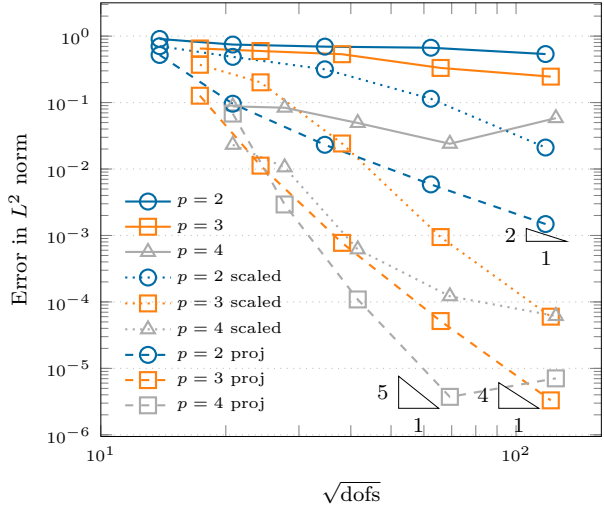
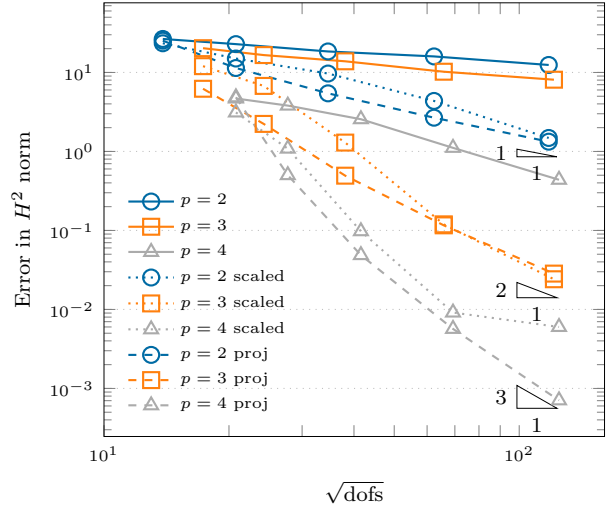
\begin{figure}
	\centering
	\begin{subfigure}[t]{0.49\textwidth}
		\centering
		%
%
%

%
%

\pgfplotsset{compat=newest}
\pgfplotsset{every tick label/.append style={font=\Huge}}
\pgfplotsset{every label/.append style={font=\Huge}}

\usetikzlibrary{calc}
\newcommand{\logLogSlopeTriangle}[5]
{

    \pgfplotsextra
    {
        \pgfkeysgetvalue{/pgfplots/xmin}{\xmin}
        \pgfkeysgetvalue{/pgfplots/xmax}{\xmax}
        \pgfkeysgetvalue{/pgfplots/ymin}{\ymin}
        \pgfkeysgetvalue{/pgfplots/ymax}{\ymax}

        \pgfmathsetmacro{\xArel}{#1}
        \pgfmathsetmacro{\yArel}{#3}
        \pgfmathsetmacro{\xBrel}{#1-#2}
        \pgfmathsetmacro{\yBrel}{\yArel}
        \pgfmathsetmacro{\xCrel}{\xArel}

        \pgfmathsetmacro{\lnxB}{\xmin*(1-(#1-#2))+\xmax*(#1-#2)} 
        \pgfmathsetmacro{\lnxA}{\xmin*(1-#1)+\xmax*#1} 
        \pgfmathsetmacro{\lnyA}{\ymin*(1-#3)+\ymax*#3} 
        \pgfmathsetmacro{\lnyC}{\lnyA+#4*(\lnxA-\lnxB)}
        \pgfmathsetmacro{\yCrel}{\lnyC-\ymin)/(\ymax-\ymin)} 

        \coordinate (A) at (rel axis cs:\xArel,\yArel);
        \coordinate (B) at (rel axis cs:\xBrel,\yBrel);
        \coordinate (C) at (rel axis cs:\xCrel,\yCrel);

        \draw[#5]   (A)-- node[pos=0.5,anchor=north] {1}
                    (B)-- 
                    (C)-- node[pos=0.5,anchor=west] {#4}
                    cycle;
    }
}

\newcommand{\logLogSlopeReverseTriangle}[5]
{
	
	\pgfplotsextra
	{
		\pgfkeysgetvalue{/pgfplots/xmin}{\xmin}
		\pgfkeysgetvalue{/pgfplots/xmax}{\xmax}
		\pgfkeysgetvalue{/pgfplots/ymin}{\ymin}
		\pgfkeysgetvalue{/pgfplots/ymax}{\ymax}
		
		\pgfmathsetmacro{\xArel}{#1}
		\pgfmathsetmacro{\yArel}{#3}
		\pgfmathsetmacro{\xBrel}{#1+#2}
		\pgfmathsetmacro{\yBrel}{\yArel}
		\pgfmathsetmacro{\xCrel}{\xArel}
		
		\pgfmathsetmacro{\lnxB}{\xmin*(1-(#1-#2))+\xmax*(#1-#2)} 
		\pgfmathsetmacro{\lnxA}{\xmin*(1-#1)+\xmax*#1} 
		\pgfmathsetmacro{\lnyA}{\ymin*(1-#3)+\ymax*#3} 
		\pgfmathsetmacro{\lnyC}{\lnyA+#4*(\lnxA-\lnxB)}
		\pgfmathsetmacro{\yCrel}{\lnyC-\ymin)/(\ymax-\ymin)}
		
		\coordinate (A) at (rel axis cs:\xArel,\yArel);
		\coordinate (B) at (rel axis cs:\xBrel,\yBrel);
		\coordinate (C) at (rel axis cs:\xCrel,\yCrel);
		
		\draw[#5]   (A)-- node[pos=0.5,anchor=north] {\scriptsize 1}
		(B)-- 
		(C)-- node[pos=0.5,anchor=east] {\scriptsize #4}
		cycle;
	}
}


\begin{tikzpicture}

\begin{axis}[%
width=0.8\textwidth,
height=5.75cm,
at={(0in,0in)},
scale only axis,
xmode=log,
xmin=10,
xtick = {10,20,30,40,50,60,70,80,90,100},
xticklabels = {$10^1$,,,,,,,,,$10^2$},
xminorticks=true,
ymode=log,
yminorticks=true,
axis background/.style={fill=white},
ylabel= Error in $L^2$ norm,
xlabel= $\sqrt{\text{dofs}}$,
legend style={font=\tiny},
legend pos= south west,
]

%
%
%
%
%

	\addplot [color=c1, solid,mark=o,mark options={solid}]
    table[ x index = 0, y index = 1] {\dataDir/shell_fourPatches_curved_breitenberger_nonconforming_p2.txt};
    \addlegendentry{$p=2$};
	
	\addplot[color=c2, solid, mark=square, mark options={solid}] 
    table[ x index = 0, y index = 1] {\dataDir/shell_fourPatches_curved_breitenberger_nonconforming_p3.txt};
    \addlegendentry{$p=3$};

	\addplot[color=c3, solid, mark=triangle, mark options={solid}] 
table[ x index = 0, y index = 1] {\dataDir/shell_fourPatches_curved_breitenberger_nonconforming_p4.txt};
\addlegendentry{$p=4$};

	\addplot [color=c1, dotted,mark=o,mark options={solid}]
    table[ x index = 0, y index = 1] {\dataDir/shell_fourPatches_curved_josef_nonconforming_p2.txt};
    \addlegendentry{$p=2$ scaled};
	
	\addplot[color=c2, dotted, mark=square, mark options={solid}] 
    table[ x index = 0, y index = 1] {\dataDir/shell_fourPatches_curved_josef_nonconforming_p3.txt};
    \addlegendentry{$p=3$ scaled};

	\addplot[color=c3, dotted, mark=triangle, mark options={solid}] 
table[ x index = 0, y index = 1] {\dataDir/shell_fourPatches_curved_josef_nonconforming_p4.txt};
\addlegendentry{$p=4$ scaled};

	\addplot [color=c1, dashed,mark=o,mark options={solid}]
    table[ x index = 0, y index = 1] {\dataDir/shell_fourPatches_curved_proj_nonconforming_p2.txt};
    \addlegendentry{$p=2$ proj};
	
	\addplot[color=c2, dashed, mark=square, mark options={solid}] 
    table[ x index = 0, y index = 1] {\dataDir/shell_fourPatches_curved_proj_nonconforming_p3.txt};
    \addlegendentry{$p=3$ proj};

	\addplot[color=c3, dashed, mark=square, mark options={solid}] 
table[ x index = 0, y index = 1] {\dataDir/shell_fourPatches_curved_proj_nonconforming_p4.txt};
\addlegendentry{$p=4$ proj};

{
\logLogSlopeReverseTriangle{0.85}{0.08}{0.45}{2}{black};
\logLogSlopeReverseTriangle{0.795}{0.08}{0.065}{4}{black};
\logLogSlopeReverseTriangle{0.595}{0.08}{0.065}{5}{black};
}
\end{axis}
\end{tikzpicture}%
		\caption{Error $L^2$.}
	\end{subfigure}
	\hfill
	\begin{subfigure}[t]{0.49\textwidth}
		\centering
		%
%
%

%
%

\pgfplotsset{compat=newest}
\pgfplotsset{every tick label/.append style={font=\Huge}}
\pgfplotsset{every label/.append style={font=\Huge}}

\usetikzlibrary{calc}
\newcommand{\logLogSlopeTriangle}[5]
{

    \pgfplotsextra
    {
        \pgfkeysgetvalue{/pgfplots/xmin}{\xmin}
        \pgfkeysgetvalue{/pgfplots/xmax}{\xmax}
        \pgfkeysgetvalue{/pgfplots/ymin}{\ymin}
        \pgfkeysgetvalue{/pgfplots/ymax}{\ymax}

        \pgfmathsetmacro{\xArel}{#1}
        \pgfmathsetmacro{\yArel}{#3}
        \pgfmathsetmacro{\xBrel}{#1-#2}
        \pgfmathsetmacro{\yBrel}{\yArel}
        \pgfmathsetmacro{\xCrel}{\xArel}

        \pgfmathsetmacro{\lnxB}{\xmin*(1-(#1-#2))+\xmax*(#1-#2)} 
        \pgfmathsetmacro{\lnxA}{\xmin*(1-#1)+\xmax*#1} 
        \pgfmathsetmacro{\lnyA}{\ymin*(1-#3)+\ymax*#3} 
        \pgfmathsetmacro{\lnyC}{\lnyA+#4*(\lnxA-\lnxB)}
        \pgfmathsetmacro{\yCrel}{\lnyC-\ymin)/(\ymax-\ymin)} 

        \coordinate (A) at (rel axis cs:\xArel,\yArel);
        \coordinate (B) at (rel axis cs:\xBrel,\yBrel);
        \coordinate (C) at (rel axis cs:\xCrel,\yCrel);

        \draw[#5]   (A)-- node[pos=0.5,anchor=north] {1}
                    (B)-- 
                    (C)-- node[pos=0.5,anchor=west] {#4}
                    cycle;
    }
}

\newcommand{\logLogSlopeReverseTriangle}[5]
{
	
	\pgfplotsextra
	{
		\pgfkeysgetvalue{/pgfplots/xmin}{\xmin}
		\pgfkeysgetvalue{/pgfplots/xmax}{\xmax}
		\pgfkeysgetvalue{/pgfplots/ymin}{\ymin}
		\pgfkeysgetvalue{/pgfplots/ymax}{\ymax}
		
		\pgfmathsetmacro{\xArel}{#1}
		\pgfmathsetmacro{\yArel}{#3}
		\pgfmathsetmacro{\xBrel}{#1+#2}
		\pgfmathsetmacro{\yBrel}{\yArel}
		\pgfmathsetmacro{\xCrel}{\xArel}
		
		\pgfmathsetmacro{\lnxB}{\xmin*(1-(#1-#2))+\xmax*(#1-#2)} 
		\pgfmathsetmacro{\lnxA}{\xmin*(1-#1)+\xmax*#1} 
		\pgfmathsetmacro{\lnyA}{\ymin*(1-#3)+\ymax*#3} 
		\pgfmathsetmacro{\lnyC}{\lnyA+#4*(\lnxA-\lnxB)}
		\pgfmathsetmacro{\yCrel}{\lnyC-\ymin)/(\ymax-\ymin)}
		
		\coordinate (A) at (rel axis cs:\xArel,\yArel);
		\coordinate (B) at (rel axis cs:\xBrel,\yBrel);
		\coordinate (C) at (rel axis cs:\xCrel,\yCrel);
		
		\draw[#5]   (A)-- node[pos=0.5,anchor=north] {\scriptsize 1}
		(B)-- 
		(C)-- node[pos=0.5,anchor=east] {\scriptsize #4}
		cycle;
	}
}


\begin{tikzpicture}

\begin{axis}[%
width=0.8\textwidth,
height=5.75cm,
at={(0in,0in)},
scale only axis,
xmode=log,
xmin=10,
xtick = {10,20,30,40,50,60,70,80,90,100},
xticklabels = {$10^1$,,,,,,,,,$10^2$},
xminorticks=true,
ymode=log,
yminorticks=true,
axis background/.style={fill=white},
ylabel= Error in $H^2$ norm,
xlabel= $\sqrt{\text{dofs}}$,
legend style={font=\tiny},
legend pos= south west,
]

	\addplot [color=c1, solid,mark=o,mark options={solid}]
    table[ x index = 0, y index = 3] {\dataDir/shell_fourPatches_curved_breitenberger_nonconforming_p2.txt};
    \addlegendentry{$p=2$};
	
	\addplot[color=c2, solid, mark=square, mark options={solid}] 
    table[ x index = 0, y index = 3] {\dataDir/shell_fourPatches_curved_breitenberger_nonconforming_p3.txt};
    \addlegendentry{$p=3$};

	\addplot[color=c3, solid, mark=triangle, mark options={solid}] 
table[ x index = 0, y index = 3] {\dataDir/shell_fourPatches_curved_breitenberger_nonconforming_p4.txt};
\addlegendentry{$p=4$};

	\addplot [color=c1, dotted,mark=o,mark options={solid}]
    table[ x index = 0, y index = 3] {\dataDir/shell_fourPatches_curved_josef_nonconforming_p2.txt};
    \addlegendentry{$p=2$ scaled};
	
	\addplot[color=c2, dotted, mark=square, mark options={solid}] 
    table[ x index = 0, y index = 3] {\dataDir/shell_fourPatches_curved_josef_nonconforming_p3.txt};
    \addlegendentry{$p=3$ scaled};

	\addplot[color=c3, dotted, mark=triangle, mark options={solid}] 
table[ x index = 0, y index = 3] {\dataDir/shell_fourPatches_curved_josef_nonconforming_p4.txt};
\addlegendentry{$p=4$ scaled};

	\addplot [color=c1, dashed,mark=o,mark options={solid}]
    table[ x index = 0, y index = 3] {\dataDir/shell_fourPatches_curved_proj_nonconforming_p2.txt};
    \addlegendentry{$p=2$ proj};
	
	\addplot[color=c2, dashed, mark=square, mark options={solid}] 
    table[ x index = 0, y index = 3] {\dataDir/shell_fourPatches_curved_proj_nonconforming_p3.txt};
    \addlegendentry{$p=3$ proj};

	\addplot[color=c3, dashed, mark=triangle, mark options={solid}] 
table[ x index = 0, y index = 3] {\dataDir/shell_fourPatches_curved_proj_nonconforming_p4.txt};
\addlegendentry{$p=4$ proj};

{
\logLogSlopeReverseTriangle{0.825}{0.08}{0.645}{1}{black};
\logLogSlopeReverseTriangle{0.825}{0.08}{0.32}{2}{black};
\logLogSlopeReverseTriangle{0.825}{0.08}{0.065}{3}{black};
}
\end{axis}
\end{tikzpicture}%
		\caption{Error $H^2$.}
	\end{subfigure}		
	\caption{Convergence study of the error measured in the $L^2$ and $H^2$ norms in the non-trimmed, non-matching four patches example for different B-splines of degree $p=2,3,4$. Comparison of a classic penalty method, the scaled version with respect to the problem parameters proposed in~\citep{Herrema2019} (\textit{scaled}) and our projection approach (\textit{proj}).}
	\label{fig:convergence_shell_curved_4patches}
\end{figure}

\noindent The results are summarized in~\Cref{fig:convergence_shell_curved_4patches}, where the convergence of the error measured in the $L^2$ and $H^2$ norms, respectively, is plotted against the square root of the number of dofs. We observe that our proposed method attains the expected order of convergence starting from very coarse meshes, whereas interface locking hinders the convergence rates of other penalty methods in the pre-asymptotic regime. As a consequence, we observe a substantial gain of accuracy per degree-of-freedom of the projection strategy, particularly in the $L^2$ norm. Additionally, we highlight the suboptimal convergence rates achieved by the method proposed in~\citep{Herrema2019}, noticeable in the asymptotic regime for $p=4$. Further, also for $p=4$ and the $L^2$ norm, we observe the detrimental impact of our choice of penalty parameters on the conditioning of the stiffness matrix and, consequently, on the solution accuracy. This effect can be mitigated by reducing the exponent $\beta$ in~\Cref{eq:projectedTerms}, knowing that the method will converge sub-optimally, as depicted in~\Cref{fig:convergence_shell_curved_4patches_exp}. For this reason, we will focus solely on moderate spline degrees $p=2,3$ in the following numerical experiments.
\begin{figure}
	\centering
	\begin{subfigure}[t]{0.49\textwidth}
		\centering
		%
%
%

%
%

\pgfplotsset{compat=newest}
\pgfplotsset{every tick label/.append style={font=\Huge}}
\pgfplotsset{every label/.append style={font=\Huge}}

\usetikzlibrary{calc}
\newcommand{\logLogSlopeTriangle}[5]
{

    \pgfplotsextra
    {
        \pgfkeysgetvalue{/pgfplots/xmin}{\xmin}
        \pgfkeysgetvalue{/pgfplots/xmax}{\xmax}
        \pgfkeysgetvalue{/pgfplots/ymin}{\ymin}
        \pgfkeysgetvalue{/pgfplots/ymax}{\ymax}

        \pgfmathsetmacro{\xArel}{#1}
        \pgfmathsetmacro{\yArel}{#3}
        \pgfmathsetmacro{\xBrel}{#1-#2}
        \pgfmathsetmacro{\yBrel}{\yArel}
        \pgfmathsetmacro{\xCrel}{\xArel}

        \pgfmathsetmacro{\lnxB}{\xmin*(1-(#1-#2))+\xmax*(#1-#2)} 
        \pgfmathsetmacro{\lnxA}{\xmin*(1-#1)+\xmax*#1} 
        \pgfmathsetmacro{\lnyA}{\ymin*(1-#3)+\ymax*#3} 
        \pgfmathsetmacro{\lnyC}{\lnyA+#4*(\lnxA-\lnxB)}
        \pgfmathsetmacro{\yCrel}{\lnyC-\ymin)/(\ymax-\ymin)} 

        \coordinate (A) at (rel axis cs:\xArel,\yArel);
        \coordinate (B) at (rel axis cs:\xBrel,\yBrel);
        \coordinate (C) at (rel axis cs:\xCrel,\yCrel);

        \draw[#5]   (A)-- node[pos=0.5,anchor=north] {1}
                    (B)-- 
                    (C)-- node[pos=0.5,anchor=west] {#4}
                    cycle;
    }
}

\newcommand{\logLogSlopeReverseTriangle}[5]
{
	
	\pgfplotsextra
	{
		\pgfkeysgetvalue{/pgfplots/xmin}{\xmin}
		\pgfkeysgetvalue{/pgfplots/xmax}{\xmax}
		\pgfkeysgetvalue{/pgfplots/ymin}{\ymin}
		\pgfkeysgetvalue{/pgfplots/ymax}{\ymax}
		
		\pgfmathsetmacro{\xArel}{#1}
		\pgfmathsetmacro{\yArel}{#3}
		\pgfmathsetmacro{\xBrel}{#1+#2}
		\pgfmathsetmacro{\yBrel}{\yArel}
		\pgfmathsetmacro{\xCrel}{\xArel}
		
		\pgfmathsetmacro{\lnxB}{\xmin*(1-(#1-#2))+\xmax*(#1-#2)} 
		\pgfmathsetmacro{\lnxA}{\xmin*(1-#1)+\xmax*#1} 
		\pgfmathsetmacro{\lnyA}{\ymin*(1-#3)+\ymax*#3} 
		\pgfmathsetmacro{\lnyC}{\lnyA+#4*(\lnxA-\lnxB)}
		\pgfmathsetmacro{\yCrel}{\lnyC-\ymin)/(\ymax-\ymin)}
		
		\coordinate (A) at (rel axis cs:\xArel,\yArel);
		\coordinate (B) at (rel axis cs:\xBrel,\yBrel);
		\coordinate (C) at (rel axis cs:\xCrel,\yCrel);
		
		\draw[#5]   (A)-- node[pos=0.5,anchor=north] {\scriptsize 1}
		(B)-- 
		(C)-- node[pos=0.5,anchor=east] {\scriptsize #4}
		cycle;
	}
}


\begin{tikzpicture}

\begin{axis}[%
width=0.8\textwidth,
height=5.75cm,
at={(0in,0in)},
scale only axis,
xmode=log,
xmin=10,
xtick = {10,20,30,40,50,60,70,80,90,100},
xticklabels = {$10^1$,,,,,,,,,$10^2$},
xminorticks=true,
ymode=log,
yminorticks=true,
axis background/.style={fill=white},
ylabel= Error in $L^2$ norm,
xlabel= $\sqrt{\text{dofs}}$,
legend style={font=\tiny},
legend style={at={(0.375,0.575)}}
]

%
%
%
%
%

	\addplot [color=c1, solid,mark=o,mark options={solid}]
    table[ x index = 0, y index = 1] {\dataDir/shell_fourPatches_curved_proj_nonconforming_p2_expPmin1.txt};
    \addlegendentry{$p=2, \beta = 1$};
	
	\addplot[color=c2, solid, mark=square, mark options={solid}] 
    table[ x index = 0, y index = 1] {\dataDir/shell_fourPatches_curved_proj_nonconforming_p3_expPmin1.txt};
    \addlegendentry{$p=3, \beta = 2$};

	\addplot[color=c3, solid, mark=triangle, mark options={solid}] 
table[ x index = 0, y index = 1] {\dataDir/shell_fourPatches_curved_proj_nonconforming_p4_expPmin1.txt};
\addlegendentry{$p=4, \beta = 3$};

	\addplot [color=c1, dotted,mark=o,mark options={solid}]
    table[ x index = 0, y index = 1] {\dataDir/shell_fourPatches_curved_proj_nonconforming_p2_expP.txt};
    \addlegendentry{$p=2, \beta = 2$};
	
	\addplot[color=c2, dotted, mark=square, mark options={solid}] 
    table[ x index = 0, y index = 1] {\dataDir/shell_fourPatches_curved_proj_nonconforming_p3_expP.txt};
    \addlegendentry{$p=3, \beta = 3$};

	\addplot[color=c3, dotted, mark=triangle, mark options={solid}] 
table[ x index = 0, y index = 1] {\dataDir/shell_fourPatches_curved_proj_nonconforming_p4_expP.txt};
\addlegendentry{$p=4, \beta = 4$};

	\addplot [color=c1, dashed,mark=o,mark options={solid}]
    table[ x index = 0, y index = 1] {\dataDir/shell_fourPatches_curved_proj_nonconforming_p2.txt};
    \addlegendentry{$p=2, \beta = 3$};
	
	\addplot[color=c2, dashed, mark=square, mark options={solid}] 
    table[ x index = 0, y index = 1] {\dataDir/shell_fourPatches_curved_proj_nonconforming_p3.txt};
    \addlegendentry{$p=3, \beta = 4$};

	\addplot[color=c3, dashed, mark=square, mark options={solid}] 
table[ x index = 0, y index = 1] {\dataDir/shell_fourPatches_curved_proj_nonconforming_p4.txt};
\addlegendentry{$p=4, \beta = 5$};

{
\logLogSlopeReverseTriangle{0.85}{0.08}{0.45}{2}{black};
\logLogSlopeReverseTriangle{0.795}{0.08}{0.065}{4}{black};
\logLogSlopeReverseTriangle{0.595}{0.08}{0.065}{5}{black};
}
\end{axis}
\end{tikzpicture}%
		\caption{Error $L^2$.}
	\end{subfigure}
	\hfill
	\begin{subfigure}[t]{0.49\textwidth}
		\centering
		%
%
%

%
%

\pgfplotsset{compat=newest}
\pgfplotsset{every tick label/.append style={font=\Huge}}
\pgfplotsset{every label/.append style={font=\Huge}}

\usetikzlibrary{calc}
\newcommand{\logLogSlopeTriangle}[5]
{

    \pgfplotsextra
    {
        \pgfkeysgetvalue{/pgfplots/xmin}{\xmin}
        \pgfkeysgetvalue{/pgfplots/xmax}{\xmax}
        \pgfkeysgetvalue{/pgfplots/ymin}{\ymin}
        \pgfkeysgetvalue{/pgfplots/ymax}{\ymax}

        \pgfmathsetmacro{\xArel}{#1}
        \pgfmathsetmacro{\yArel}{#3}
        \pgfmathsetmacro{\xBrel}{#1-#2}
        \pgfmathsetmacro{\yBrel}{\yArel}
        \pgfmathsetmacro{\xCrel}{\xArel}

        \pgfmathsetmacro{\lnxB}{\xmin*(1-(#1-#2))+\xmax*(#1-#2)} 
        \pgfmathsetmacro{\lnxA}{\xmin*(1-#1)+\xmax*#1} 
        \pgfmathsetmacro{\lnyA}{\ymin*(1-#3)+\ymax*#3} 
        \pgfmathsetmacro{\lnyC}{\lnyA+#4*(\lnxA-\lnxB)}
        \pgfmathsetmacro{\yCrel}{\lnyC-\ymin)/(\ymax-\ymin)} 

        \coordinate (A) at (rel axis cs:\xArel,\yArel);
        \coordinate (B) at (rel axis cs:\xBrel,\yBrel);
        \coordinate (C) at (rel axis cs:\xCrel,\yCrel);

        \draw[#5]   (A)-- node[pos=0.5,anchor=north] {1}
                    (B)-- 
                    (C)-- node[pos=0.5,anchor=west] {#4}
                    cycle;
    }
}

\newcommand{\logLogSlopeReverseTriangle}[5]
{
	
	\pgfplotsextra
	{
		\pgfkeysgetvalue{/pgfplots/xmin}{\xmin}
		\pgfkeysgetvalue{/pgfplots/xmax}{\xmax}
		\pgfkeysgetvalue{/pgfplots/ymin}{\ymin}
		\pgfkeysgetvalue{/pgfplots/ymax}{\ymax}
		
		\pgfmathsetmacro{\xArel}{#1}
		\pgfmathsetmacro{\yArel}{#3}
		\pgfmathsetmacro{\xBrel}{#1+#2}
		\pgfmathsetmacro{\yBrel}{\yArel}
		\pgfmathsetmacro{\xCrel}{\xArel}
		
		\pgfmathsetmacro{\lnxB}{\xmin*(1-(#1-#2))+\xmax*(#1-#2)} 
		\pgfmathsetmacro{\lnxA}{\xmin*(1-#1)+\xmax*#1} 
		\pgfmathsetmacro{\lnyA}{\ymin*(1-#3)+\ymax*#3} 
		\pgfmathsetmacro{\lnyC}{\lnyA+#4*(\lnxA-\lnxB)}
		\pgfmathsetmacro{\yCrel}{\lnyC-\ymin)/(\ymax-\ymin)}
		
		\coordinate (A) at (rel axis cs:\xArel,\yArel);
		\coordinate (B) at (rel axis cs:\xBrel,\yBrel);
		\coordinate (C) at (rel axis cs:\xCrel,\yCrel);
		
		\draw[#5]   (A)-- node[pos=0.5,anchor=north] {\scriptsize 1}
		(B)-- 
		(C)-- node[pos=0.5,anchor=east] {\scriptsize #4}
		cycle;
	}
}


\begin{tikzpicture}

\begin{axis}[%
width=0.8\textwidth,
height=5.75cm,
at={(0in,0in)},
scale only axis,
xmode=log,
xmin=10,
xtick = {10,20,30,40,50,60,70,80,90,100},
xticklabels = {$10^1$,,,,,,,,,$10^2$},
xminorticks=true,
ymode=log,
yminorticks=true,
axis background/.style={fill=white},
ylabel= Error in $H^2$ norm,
xlabel= $\sqrt{\text{dofs}}$,
legend style={font=\tiny},
legend pos= south west,
]

	\addplot [color=c1, solid,mark=o,mark options={solid}]
table[ x index = 0, y index = 3] {\dataDir/shell_fourPatches_curved_proj_nonconforming_p2_expPmin1.txt};
\addlegendentry{$p=2, \beta = 1$};

\addplot[color=c2, solid, mark=square, mark options={solid}] 
table[ x index = 0, y index = 3] {\dataDir/shell_fourPatches_curved_proj_nonconforming_p3_expPmin1.txt};
\addlegendentry{$p=3, \beta = 2$};

\addplot[color=c3, solid, mark=triangle, mark options={solid}] 
table[ x index = 0, y index = 3] {\dataDir/shell_fourPatches_curved_proj_nonconforming_p4_expPmin1.txt};
\addlegendentry{$p=4, \beta = 3$};

\addplot [color=c1, dotted,mark=o,mark options={solid}]
table[ x index = 0, y index = 3] {\dataDir/shell_fourPatches_curved_proj_nonconforming_p2_expP.txt};
\addlegendentry{$p=2, \beta = 2$};

\addplot[color=c2, dotted, mark=square, mark options={solid}] 
table[ x index = 0, y index = 3] {\dataDir/shell_fourPatches_curved_proj_nonconforming_p3_expP.txt};
\addlegendentry{$p=3, \beta = 3$};

\addplot[color=c3, dotted, mark=triangle, mark options={solid}] 
table[ x index = 0, y index = 3] {\dataDir/shell_fourPatches_curved_proj_nonconforming_p4_expP.txt};
\addlegendentry{$p=4, \beta = 4$};

\addplot [color=c1, dashed,mark=o,mark options={solid}]
table[ x index = 0, y index = 3] {\dataDir/shell_fourPatches_curved_proj_nonconforming_p2.txt};
\addlegendentry{$p=2, \beta = 3$};

\addplot[color=c2, dashed, mark=square, mark options={solid}] 
table[ x index = 0, y index = 3] {\dataDir/shell_fourPatches_curved_proj_nonconforming_p3.txt};
\addlegendentry{$p=3, \beta = 4$};

\addplot[color=c3, dashed, mark=square, mark options={solid}] 
table[ x index = 0, y index = 3] {\dataDir/shell_fourPatches_curved_proj_nonconforming_p4.txt};
\addlegendentry{$p=4, \beta = 5$};

{
\logLogSlopeReverseTriangle{0.825}{0.08}{0.645}{1}{black};
\logLogSlopeReverseTriangle{0.825}{0.08}{0.33}{2}{black};
\logLogSlopeReverseTriangle{0.80}{0.08}{0.065}{3}{black};
}
\end{axis}
\end{tikzpicture}%
		\caption{Error $H^2$.}
	\end{subfigure}		
	\caption{Convergence study of the error measured in the $L^2$ and $H^2$ norms for the projection method in the non-trimmed, non-matching four patches example for different B-splines of degree $p=2,3,4$. Comparison of various scaling exponent $\beta$ of the penalty parameters in~\Cref{eq:projectedTerms}.}
	\label{fig:convergence_shell_curved_4patches_exp}
\end{figure}
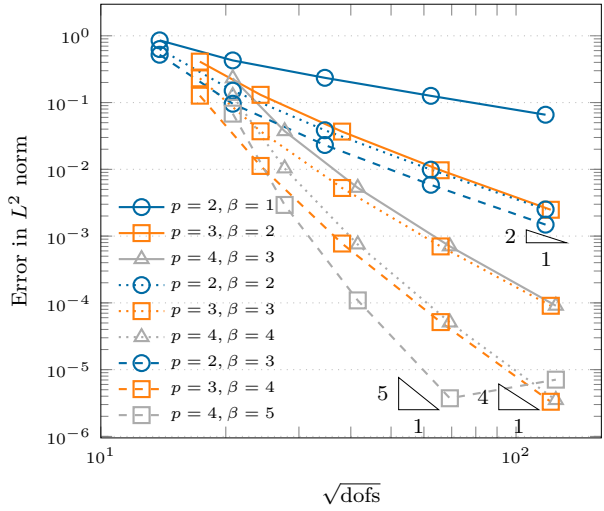
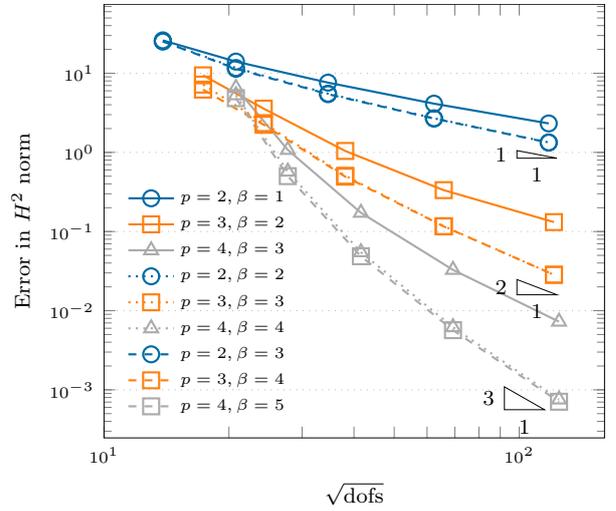
\begin{remark}
In order to retain optimal rates of convergence, whenever a cross-point is present in the geometry, we must impose a $C^0$-continuity constraint at the cross-point. For further details and a possible implementation, we refer to~\citep{Coradello2020projected}.
\end{remark}

\subsection{Scordelis-Lo roof}
\begin{figure}
	\centering
	\begin{subfigure}[t]{0.495\textwidth}
		\centering
		\includegraphics[width=\textwidth]{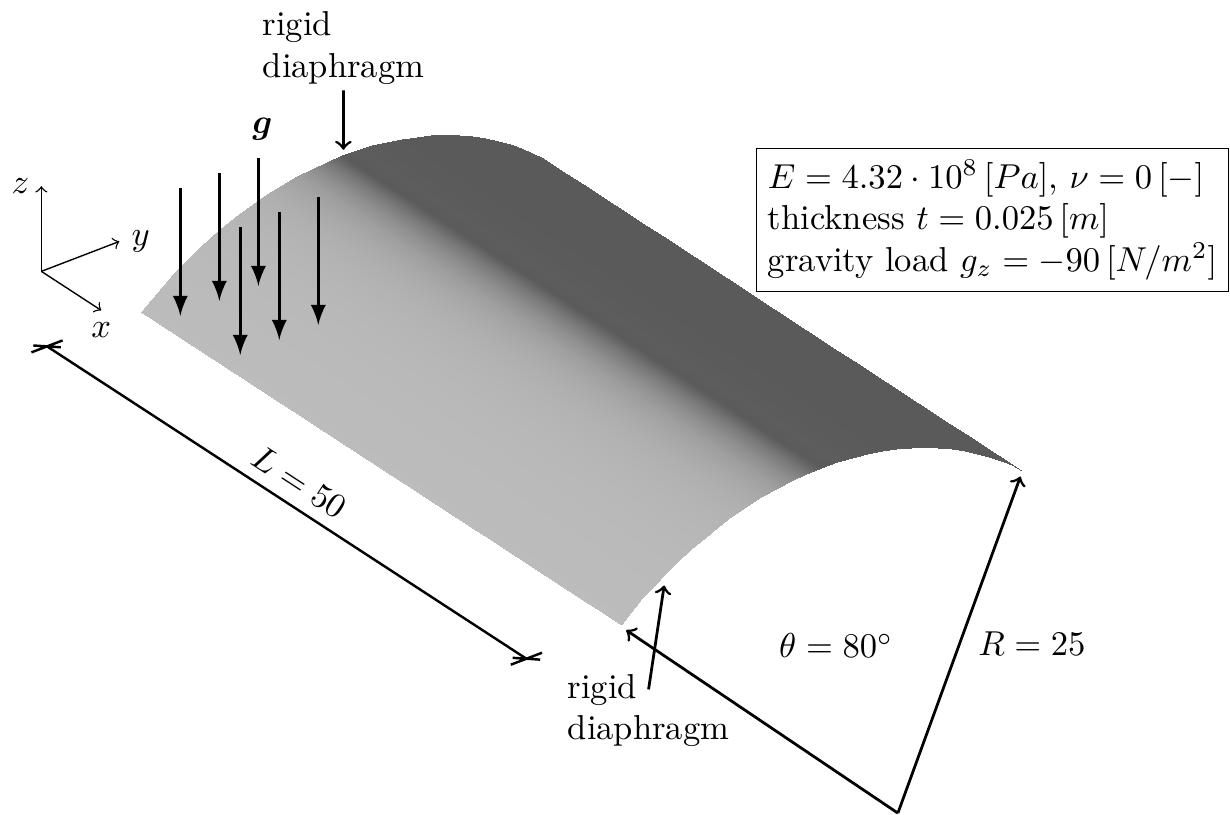}
		\caption{Geometry setup and physical parameters.}
	\end{subfigure}
	\hfill
	\begin{subfigure}[t]{0.495\textwidth}
	\centering
	\includegraphics[width=\textwidth]{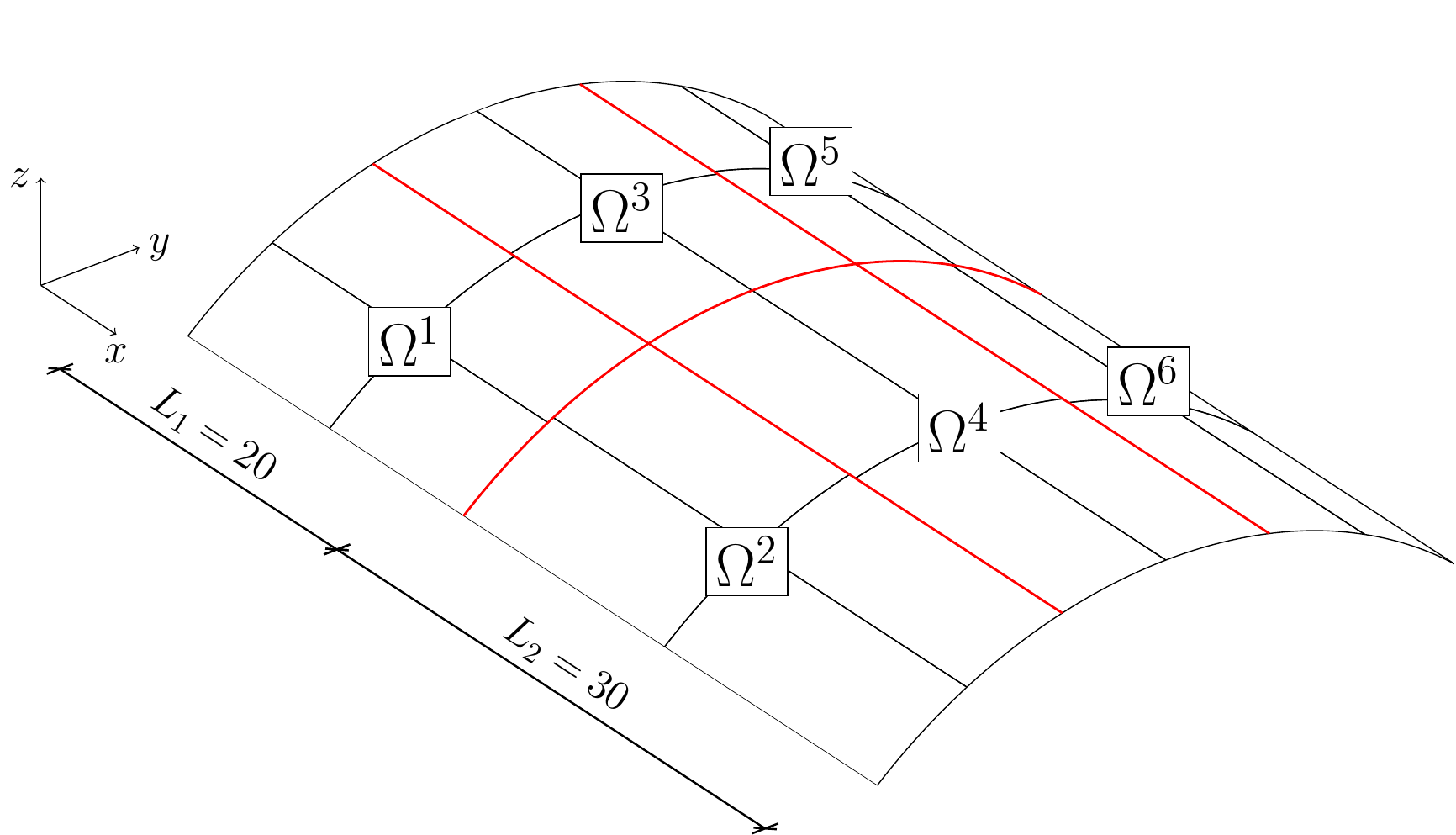}
	\caption{Initial discretization, patch interfaces are colored in red.}
\end{subfigure}
	\caption{Problem setup for the Scordelis-Lo roof example.}
\label{fig:scordelis_6patches_setup}
\end{figure}

\begin{figure}
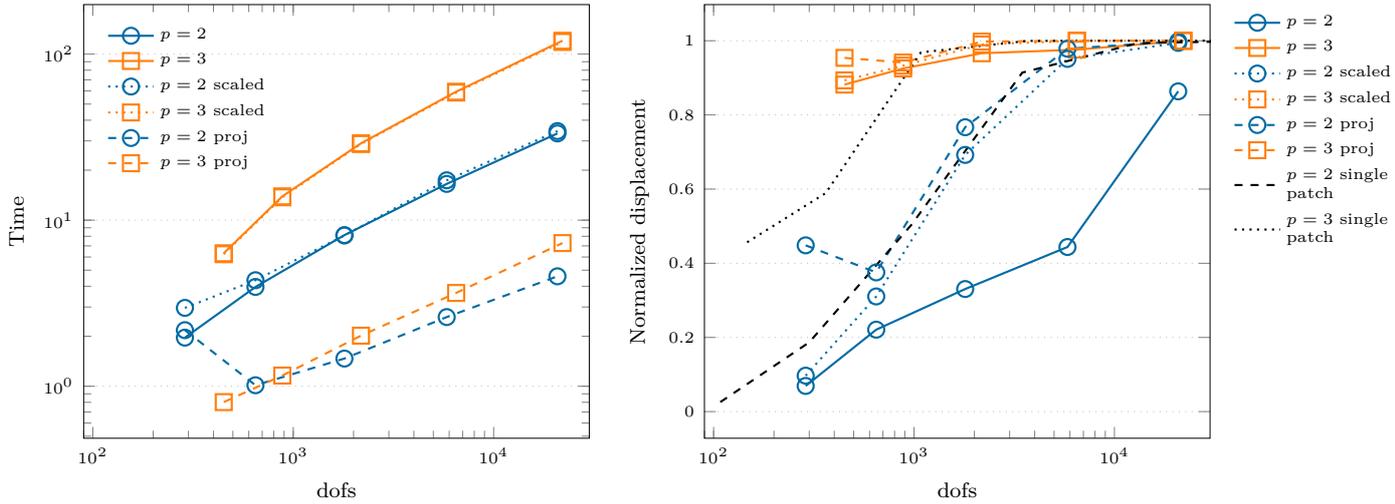

	\begin{subfigure}[t]{0.495\textwidth}
	\centering
	\input{\graphDir/runtime_scordelis_6patches.tex}
	\caption{Time needed to assemble the coupling terms.}\label{fig:runtime_scordelis_6patches}
\end{subfigure}	
	\begin{subfigure}[t]{0.495\textwidth}
		\centering
		\input{\graphDir/convergence_disp_scordelis_6patches.tex}
		\caption{Convergence of the normalized displacement in the $z$-direction.}\label{fig:scordelis_6patches_disp}
	\end{subfigure}	
	\caption{Convergence of the normalized displacement in the $z$-direction at the middle of the free edge and time needed to assemble the penalty contributions, Scordelis-Lo roof example.}
	\label{fig:scordelis_6patches}
\end{figure}

In this example we asses the performance of our method on the well-known Scordelis-Lo roof, firstly introduced as part of the shell obstacle course in~\citep{Belytschko1985}. The geometrical setup, the chosen parameters and the initial non-conforming multi-patch design, where the roof is split into six subdomains $\Omega^i, i=1,\ldots,6$,  are summarized in~\Cref{fig:scordelis_6patches_setup}. The structure is supported at both ends of the cylindrical roof by so-called rigid diaphragms, which fix the displacement in the $y$ and $z$ directions, respectively. Moreover, the roof is subjected to a uniform gravity load, directed in the negative $z$-direction.  As studied in~\citep{Herrema2019}, we modify the original thickness of the benchmark problem. In particular, we set the Young's modulus, the Poisson's ratio and the thickness of the structure to $4.32 \cdot 10^8 \, [Pa]$, $0.0 \, [-]$ and $0.025 \, [m]$, respectively. As typically done for this problem, we study the convergence of the displacement in the $z$-direction at the center of the free edge, where the reference value $u_z^{\text{ref}} = -32.01045$ is used for normalization. The results are presented in~\Cref{fig:scordelis_6patches_disp} for different penalty methods and also for the single-patch case. We observe that our approach, the method presented in~\citep{Herrema2019} and the single patch case show a similar convergence behavior. However, in case when the penalty parameter is only scaled by the Young's modulus, interface locking phenomena arise, and they are particularly severe for quadratic B-splines. Moreover, in~\Cref{fig:runtime_scordelis_6patches}, we compare the time needed to compute and assemble the penalty terms for the aforementioned approaches, where the projection method shows its computational efficiency. This is linked to the fact that, although the projection algorithm requires the solution of an additional system, the corresponding coupling terms involve significantly fewer dofs compared to standard penalty-like methods.
\begin{remark}
Although the Scordelis-Lo roof is a classical benchmark for shell analysis, it only provides a reference value for the displacement in a point. Therefore, it is not suited to quantify the order of convergence of a method, but it only serves as verification of the latter. 
\end{remark}

\subsection{L-beam}
\begin{figure}
	\centering
	\begin{subfigure}[t]{0.495\textwidth}
		\centering
		\includegraphics[width=\textwidth]{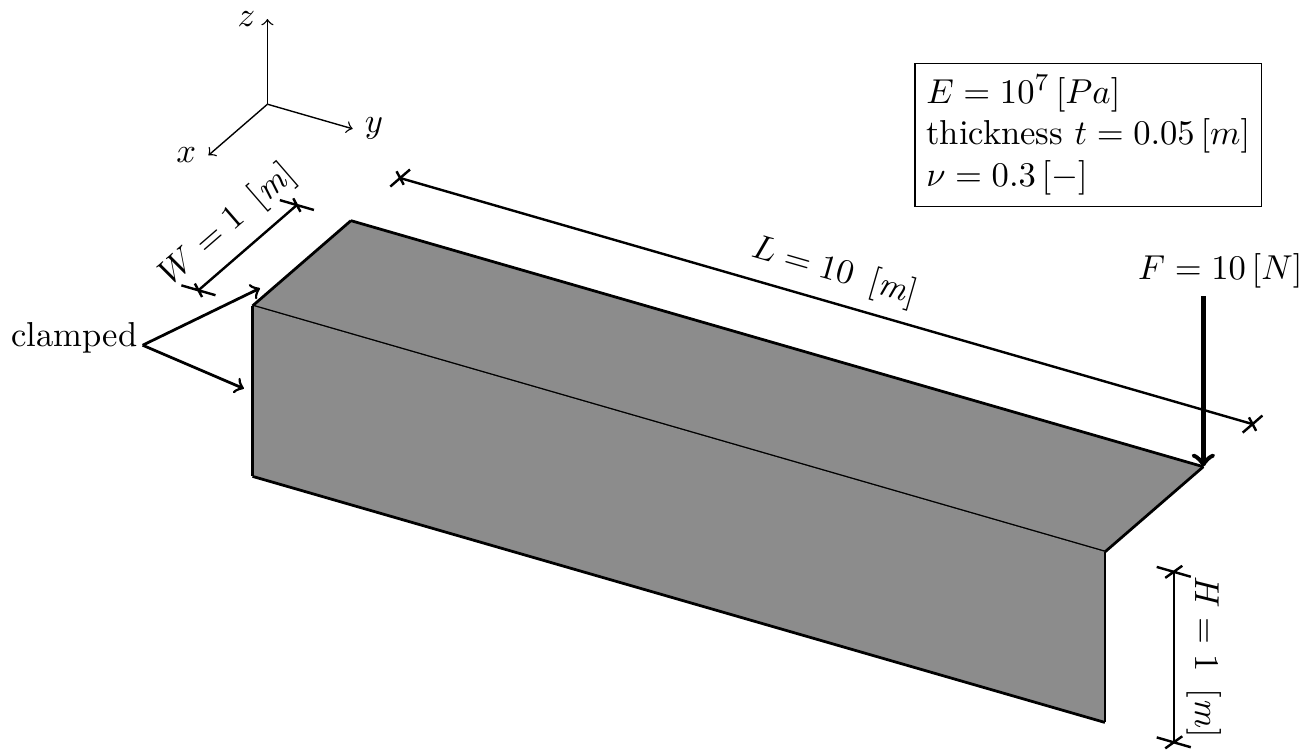}
		\caption{Geometry setup and physical parameters.}
	\end{subfigure}
	\hfill
	\begin{subfigure}[t]{0.495\textwidth}
		\centering
		\includegraphics[width=0.8\textwidth]{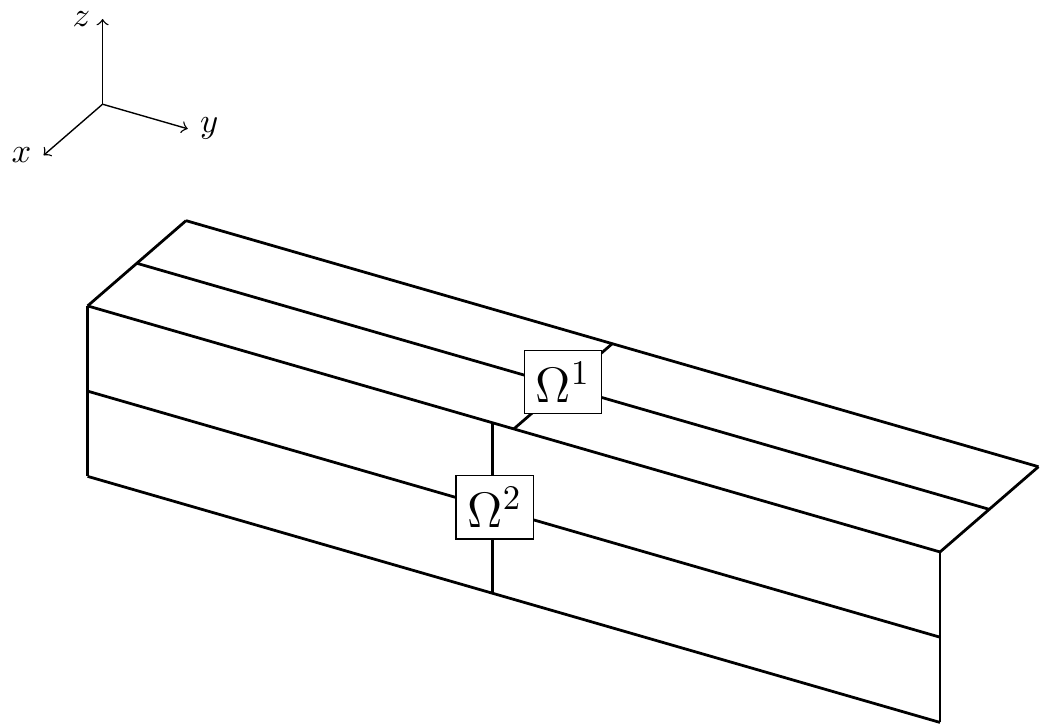}
		\caption{Initial discretization.}
	\end{subfigure}	
	\begin{subfigure}[t]{0.55\textwidth}
	\centering
	\includegraphics[width=\textwidth]{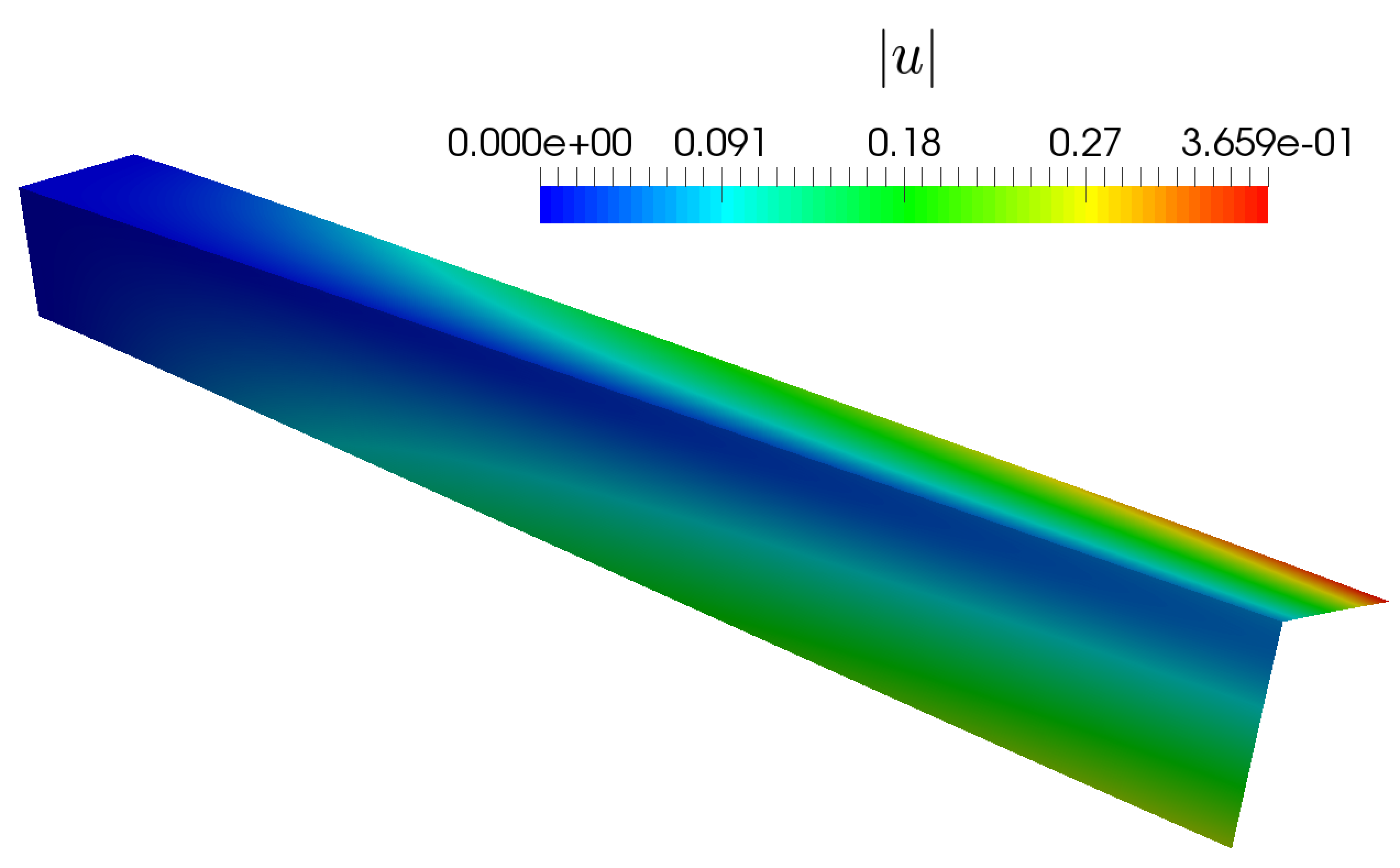}
	\caption{$L^2$ norm of the displacement, warped by a factor of 1. \mbox{B-splines} of degree $p=2$, 8192 elements.}
\end{subfigure}	
	\caption{Problem setup for the L-beam example and $L^2$ norm of the displacement field obtained with the proposed coupling strategy.}
	\label{fig:setup_l_beam}
\end{figure}

\begin{figure}
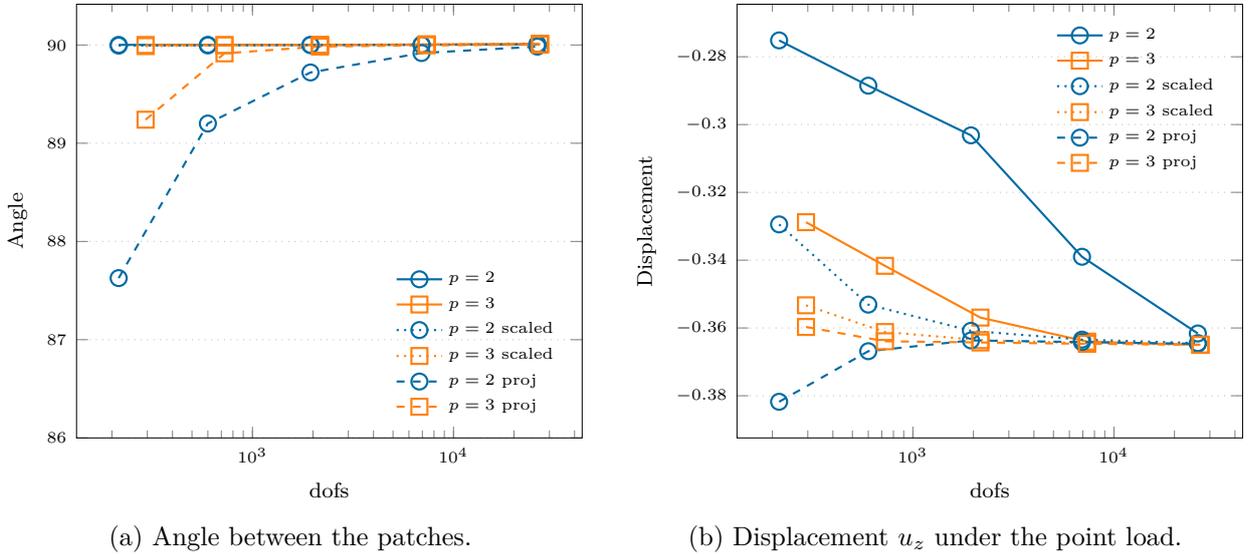

	\centering
	\begin{subfigure}[t]{0.495\textwidth}
		\centering
		\input{\graphDir/convergence_angle_l_beam.tex}
		\caption{Angle between the patches.}\label{fig:angle_l_beam}
	\end{subfigure}
	\hfill
	\begin{subfigure}[t]{0.495\textwidth}
		\centering
		\input{\graphDir/convergence_disp_l_beam.tex}
		\caption{Displacement $u_z$ under the point load.}\label{fig:disp_l_beam}
	\end{subfigure}	
	\caption{Convergence study of the angle between the patches and the displacement $u_z$ under the point load in the L-beam example, B-splines of degree $p=2,3$. Comparison of a classic penalty method, the scaled version with respect to the problem parameters proposed in~\citep{Herrema2019} (\textit{scaled}) and our projection approach (\textit{proj}).}
	\label{fig:convergence_l_beam}
\end{figure}

This example is meant to demonstrate the applicability of our approach to couple patches at an arbitrary angle, where the corresponding rotational constraint keeps the angle fix during deformation. We consider a beam with an L-section discretized by two non-conforming patches $\Omega^i, i=1,2$, as depicted in~\Cref{fig:setup_l_beam}. The beam is clamped on one side and it is subjected to a point load of $10 \, [N]$, directed in the negative $z$-direction. Further, we set the Young's modulus, the Poisson's ratio and the thickness of the structure to $10^7 \, [Pa]$, $0.3 \, [-]$ and $0.05 \, [m]$, respectively.
To check the correct imposition of the rotational constraint, we compute the angle formed by the two patches at the free corner on a series of uniformly refined meshes. The corresponding results are presented in~\Cref{fig:angle_l_beam}. We remark that on coarse meshes and for the projection method, the rotational constraint is imposed in a less ``rigid'' way compared to other penalty approaches. This allows to mitigate the effects related to interface locking starting from coarse meshes. 
Similarly to the previous example, we observe a faster convergence behavior of the vertical displacement under the point load when our approach is employed, see~\Cref{fig:disp_l_beam}, especially compared to a classical penalty method.

\subsection{Pure bending of three trimmed planar patches}
\begin{figure}
	\centering
\begin{subfigure}[t]{0.495\textwidth}
	\centering
	\includegraphics[width=0.9\textwidth]{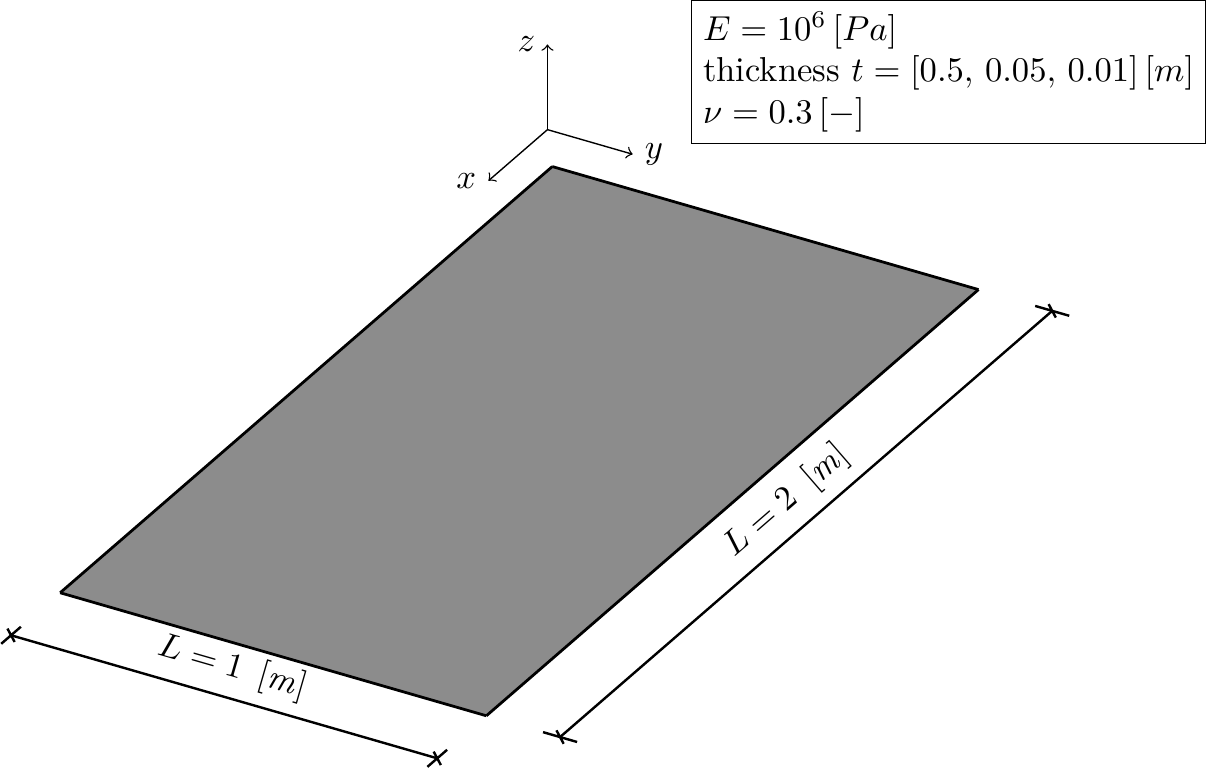}
	\caption{Geometry setup and physical parameters.}
\end{subfigure}
\hfill
\begin{subfigure}[t]{0.495\textwidth}
	\centering
	\includegraphics[width=\textwidth]{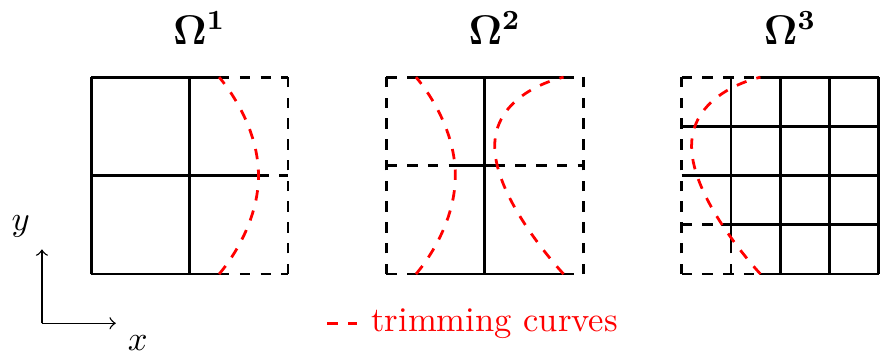}
	\caption{Initial discretization.}
\end{subfigure}	
\caption{Problem setup for the trimmed three planar patches example. For the sake of visualization, the patches have been separated.}
\label{fig:setup_trimmed_curved_3patches}
\end{figure}

\begin{figure}
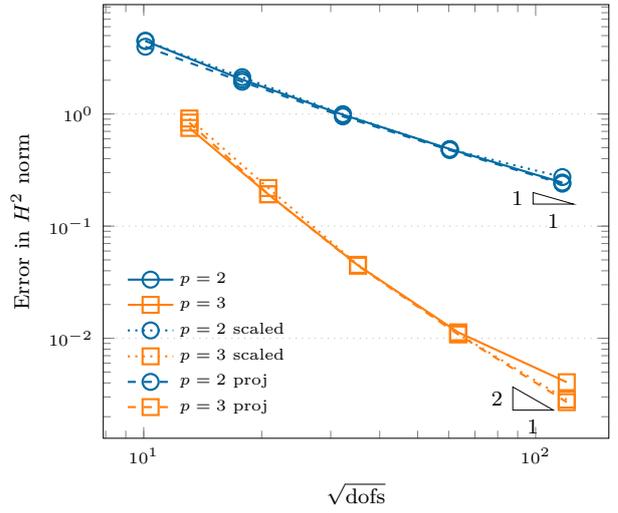
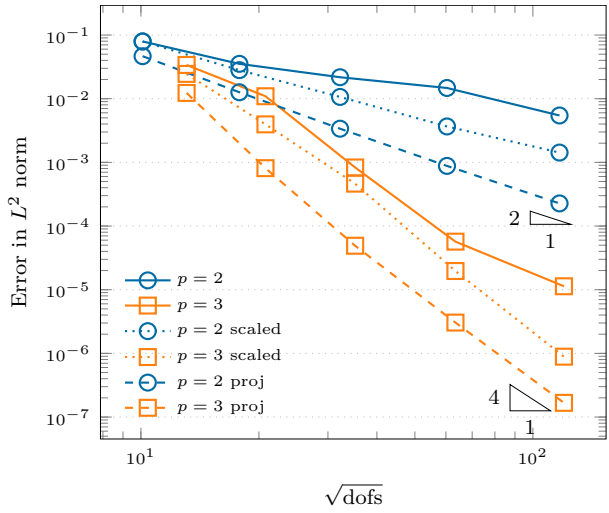
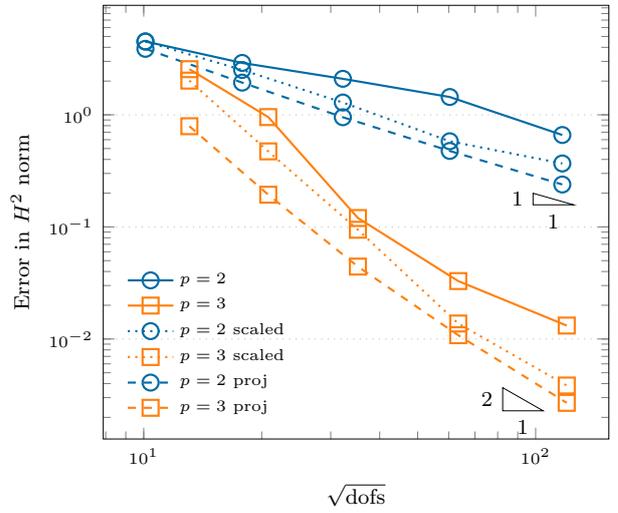
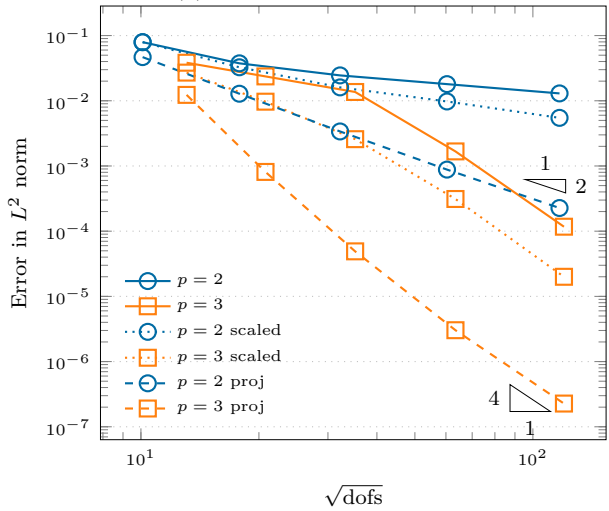
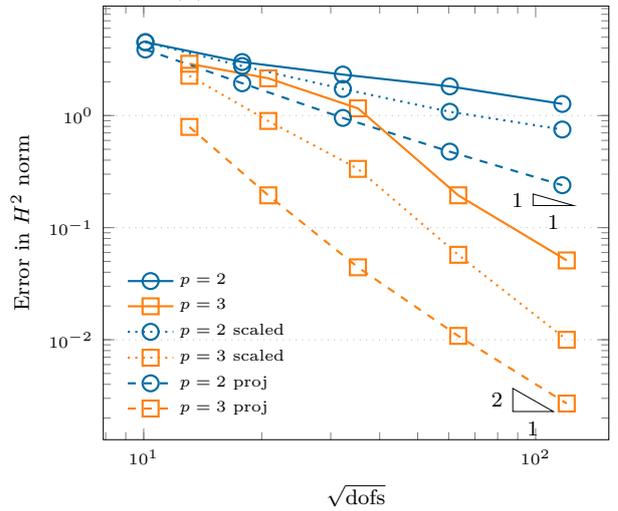

	\centering
	\begin{subfigure}[t]{0.495\textwidth}
		\centering
		\input{\graphDir/convergence_3plates_trimmed_nonconforming_L2_t5.tex}
		\caption{$L^2$ norm, $t=0.5$.}
	\end{subfigure}
	\hfill
	\begin{subfigure}[t]{0.495\textwidth}
		\centering
		\input{\graphDir/convergence_3plates_trimmed_nonconforming_H2_t5.tex}
		\caption{$H^2$ norm, $t=0.5$.}
	\end{subfigure}	
	\begin{subfigure}[t]{0.495\textwidth}
	\centering
	\input{\graphDir/convergence_3plates_trimmed_nonconforming_L2_t05.tex}
	\caption{$L^2$ norm, $t=0.05$.}
	\end{subfigure}
	\hfill
	\begin{subfigure}[t]{0.495\textwidth}
		\centering
		\input{\graphDir/convergence_3plates_trimmed_nonconforming_H2_t05.tex}
		\caption{$H^2$ norm, $t=0.05$.}
	\end{subfigure}	
	\begin{subfigure}[t]{0.495\textwidth}
	\centering
	\input{\graphDir/convergence_3plates_trimmed_nonconforming_L2_t005.tex}
	\caption{$L^2$ norm, $t=0.01$.}
	\end{subfigure}
	\hfill
	\begin{subfigure}[t]{0.495\textwidth}
		\centering
		\input{\graphDir/convergence_3plates_trimmed_nonconforming_H2_t005.tex}
		\caption{$H^2$ norm, $t=0.01$.}
	\end{subfigure}	
	\caption{Convergence study of the error measured in the $H^2$ and $L^2$ norms in the trimmed three patches example, B-splines of degree $p=2,3$, thickness $t=[0.5 \, , 0.05 \, , 0.01]  \, [m]$. Comparison of a classic penalty method, the scaled version with respect to the problem parameters proposed in~\citep{Herrema2019} (\textit{scaled}) and our projection approach (\textit{proj}).}
	\label{fig:convergence_3plates_nonconforming}
\end{figure}

In this example we consider the computational domain $\Omega = [0,2] \times [0,1]$ split into three trimmed subdomains $\Omega^i, i=1,2,3$ as depicted in~\Cref{fig:setup_trimmed_curved_3patches}. 
We remark that in this particular setup, the middle patch is coupled on both sides along trimming interfaces, defined by quadratic spline curves. 
The applied boundary conditions and loading function are again derived from a smooth solution of the form:
\begin{align}
	\bm{u}^{\text{ex}} (x, y, z) = 
	\begin{pmatrix}
		u^x \\
		u^y \\
		u^z 
	\end{pmatrix} = 
	\begin{pmatrix}
		0 \\
		0\\
		\sin(\pi  x) \sin(\pi  y) 
	\end{pmatrix} 
	\, .
\end{align} 
Then, we fix the Young's modulus and the Poisson's ratio of the structure to $10^6 \, [Pa]$ and $0.3 \, [-]$, respectively.
This example confirms the severity of locking interface phenomena, especially as the shell gets progressively slender. Indeed, we vary the thickness in the range $[0.5 \, , 0.05 \, , 0.01]  \, [m]$, where the results are reported in~\Cref{fig:convergence_3plates_nonconforming}. Moreover, we observe that in the trimmed case these detrimental effects are even more pronounced, since more basis functions are involved in the imposition of the constraints compared to the non-trimmed case. 
\begin{remark}
In the trimmed case, at any point of a coupling interface $\gamma^\ell$ and for each neighboring patch, we have $(p+1) \times (p+1)$ shape functions providing a non-zero contribution to the penalty matrices. This is in contrast to the non-trimmed case, where at any point of $\gamma^\ell$ and for each neighboring patch we have at most $(p+1)$, respectively, $2 (p+1)$ B-splines involved in the computation of the displacement and rotational coupling terms. 
\end{remark}

\subsection{Trimmed astroid}

\begin{figure}
	\centering
	\begin{subfigure}[t]{0.45\textwidth}
		\centering
	\includegraphics[width=0.8\textwidth]{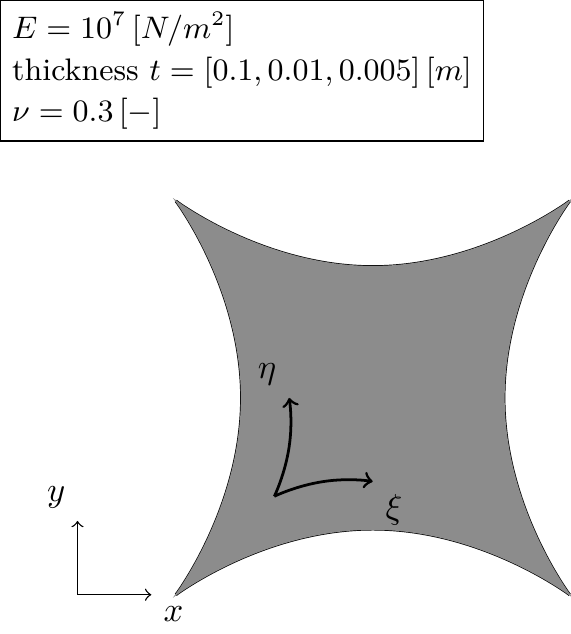}
		\caption{Geometry setup and physical parameters.}
	\end{subfigure}
	\hfill
	\begin{subfigure}[t]{0.495\textwidth}
		\centering
	\includegraphics[width=0.7\textwidth]{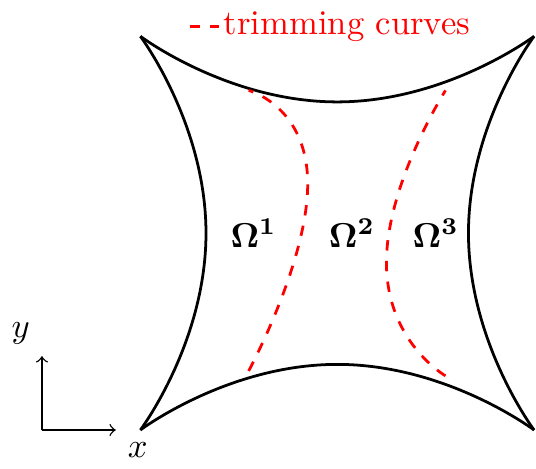}
		\caption{Trimmed multi-patch geometry.}
	\end{subfigure}	
	\caption{Problem setup and initial multi-patch geometry for the trimmed astroid example.}
\label{fig:setup_astroid}
\end{figure}

\begin{figure}
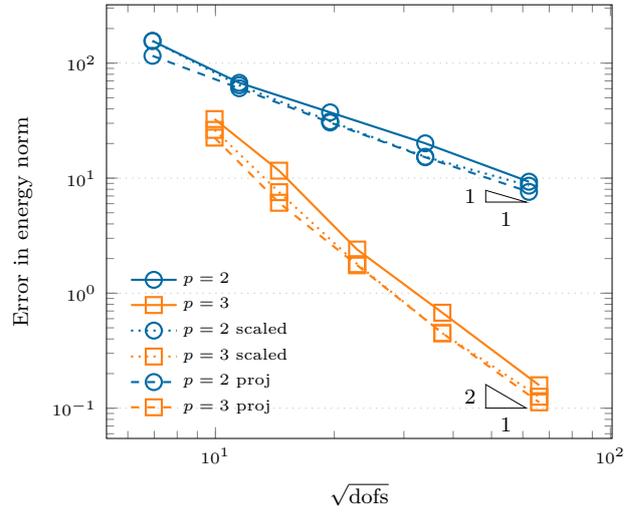
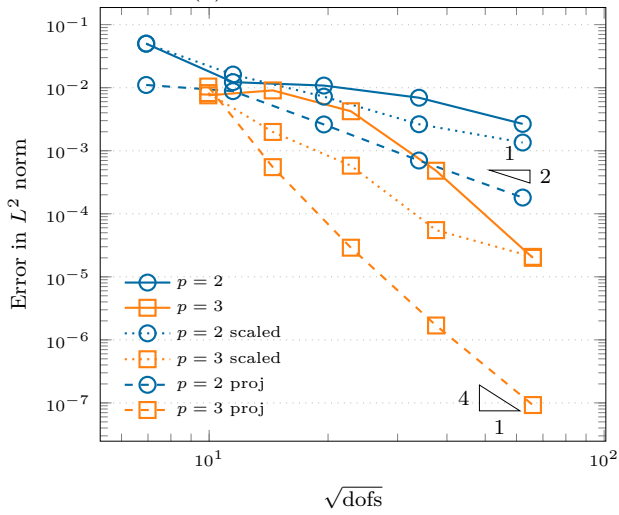
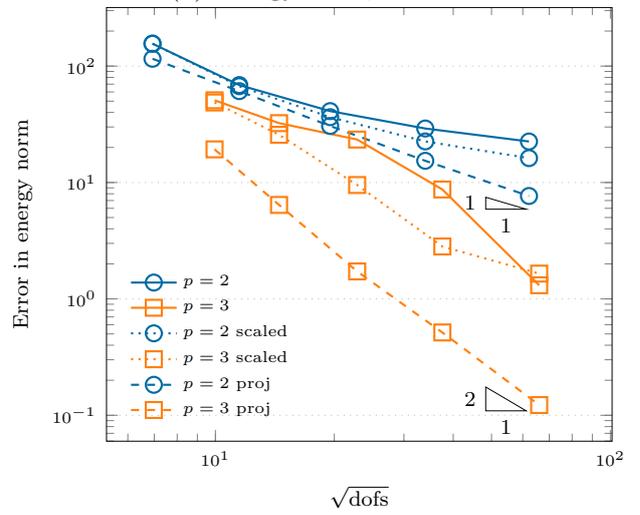
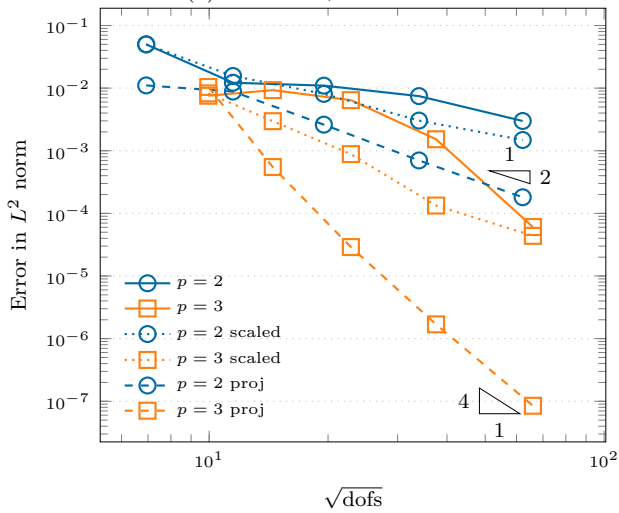
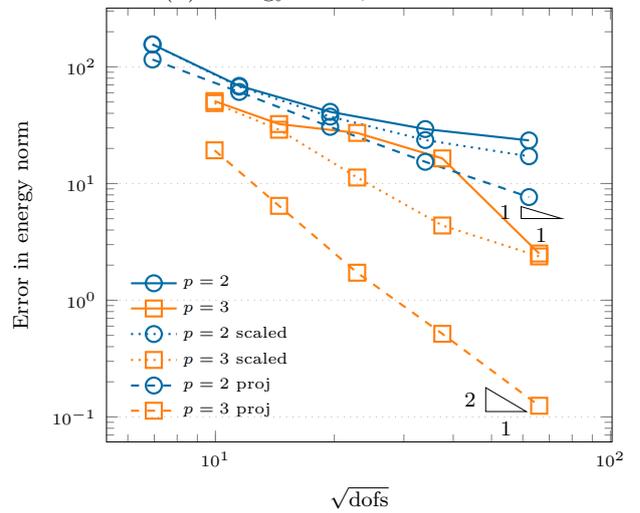

	\centering
	\begin{subfigure}[t]{0.495\textwidth}
		\centering
		\input{\graphDir/convergence_astroid_curved_nonconforming_L2_t1.tex}
		\caption{$L^2$ norm, $t=0.1$.}
	\end{subfigure}
	\hfill
	\begin{subfigure}[t]{0.495\textwidth}
		\centering
		\input{\graphDir/convergence_astroid_curved_nonconforming_energy_t1.tex}
		\caption{Energy norm, $t=0.1$.}
	\end{subfigure}	
	\begin{subfigure}[t]{0.495\textwidth}
		\centering
		\input{\graphDir/convergence_astroid_curved_nonconforming_L2_t01.tex}
		\caption{$L^2$ norm, $t=0.01$.}
	\end{subfigure}
	\hfill
	\begin{subfigure}[t]{0.495\textwidth}
		\centering
		\input{\graphDir/convergence_astroid_curved_nonconforming_energy_t01.tex}
		\caption{Energy norm, $t=0.01$.}
	\end{subfigure}	
	\begin{subfigure}[t]{0.495\textwidth}
		\centering
		\input{\graphDir/convergence_astroid_curved_nonconforming_L2_t005.tex}
		\caption{$L^2$ norm, $t=0.005$.}
	\end{subfigure}
	\hfill
	\begin{subfigure}[t]{0.495\textwidth}
		\centering
		\input{\graphDir/convergence_astroid_curved_nonconforming_energy_t005.tex}
		\caption{Energy norm, $t=0.005$.}
	\end{subfigure}	
	\caption{Convergence study of the error measured in the energy and $L^2$ norms in the trimmed astroid example, B-splines of degree $p=2,3$, thickness $t=[0.1 \, , 0.01 \, , 0.005] \, [m]$. Comparison of a classic penalty method, the scaled version with respect to the problem parameters proposed in~\citep{Herrema2019} (\textit{scaled}) and our projection approach (\textit{proj}).}
	\label{fig:convergence_astroid}
\end{figure}

This example is adapted from the shell obstacle course presented in~\citep{Benzaken2020}. We consider the computational domain $\Omega$ split into three trimmed subdomains $\Omega^i, i=1,2,3$ as depicted in~\Cref{fig:setup_astroid}. In the same figure, the two trimmed interfaces $\gamma^\ell, \ell=1,2$ are defined as quadratic B-spline curves. 
The domain $\Omega$ is characterized by the control points $\mathbf{P}_{ij}, \, i,j=1,\ldots,3$ as summarized in~\Cref{tab:astroid}, where the indices $i,j$ are ordered as the parametric coordinates $\xi, \eta$ represented in~\Cref{fig:setup_astroid}.
\begin{table}
	\centering
\begin{tabular}{cccccccccc}
	\hline $\mathbf{P}_{ij}$ & 11 & 12 & 13 & 21 & 22 & 23 & 31 & 32 & 33 \\
	\hline $x$ & 0 & 1/3 & 0 & 1/2 & 1/2 & 1/2 & 1 & 2/3 & 1\\ 
	$y$ & 0 & 1/2 & 1 & 1/3 & 1/2 & 2/3 & 0 & 1/2 & 1\\ 
	$z$ & 0 & 0 & 0 & 0 & 0 & 0 & 0 & 0 & 0\\ 
	\hline
\end{tabular}
\caption{Coordinates of the control points $\mathbf{P}_{ij}, \, i,j=1,\ldots,3$ associated to the astroid domain.}\label{tab:astroid}
\end{table}
The load function and boundary data are computed from the following manufactured solution:
\begin{align}
\bm{u}^{\text{ex}} (\xi, \eta) = 
\begin{pmatrix}
u^x \\
u^y \\
u^z 
\end{pmatrix} = 
\begin{pmatrix}
	(\frac{1}{2} - \eta)  \xi^2  (\xi - 1)^2  \eta  (1 - \eta) \\
	(\xi - \frac{1}{2})  \eta^2  (\eta - 1)^2  \xi  (1 - \xi)\\
	\xi (1 - \xi) \sin(\pi  \xi) \sin(\pi  \eta) 
\end{pmatrix} 
 \, ,
\end{align}
where we impose inhomogeneous Dirichlet and Neumann type boundary conditions on the displacements and on the bending moments, respectively, on the entire boundary $\partial \Omega$. Even though this problem is defined on a planar geometry, meaning that bending and membrane responses are decoupled, its investigation is still worthwhile since the solution $\bm{u}$ is defined as a function of the parametrization. This drastically complicates the derivation of the exact quantities and their stable computation.
%
%
\noindent The convergence results for the error in the $L^2$ and energy norms for several values of the thickness $t=\left[0.1\, , 0.01 \, , 0.005\right]  \, [m]$ are depicted in~\Cref{fig:convergence_astroid}. This example confirms that our projection method mitigates the detrimental effects linked to interface locking, yielding a significant gain of accuracy per-degrees-of-freedom. This is particularly noticeable as the thickness of the structure becomes smaller, where, for other penalty techniques, locking phenomena hinder the optimal convergence in the pre-asymptotic regime.

\subsection{Trimmed cylinder}
\begin{figure}
	\centering
	\begin{subfigure}[t]{0.495\textwidth}
		\centering
		\includegraphics[width=\textwidth]{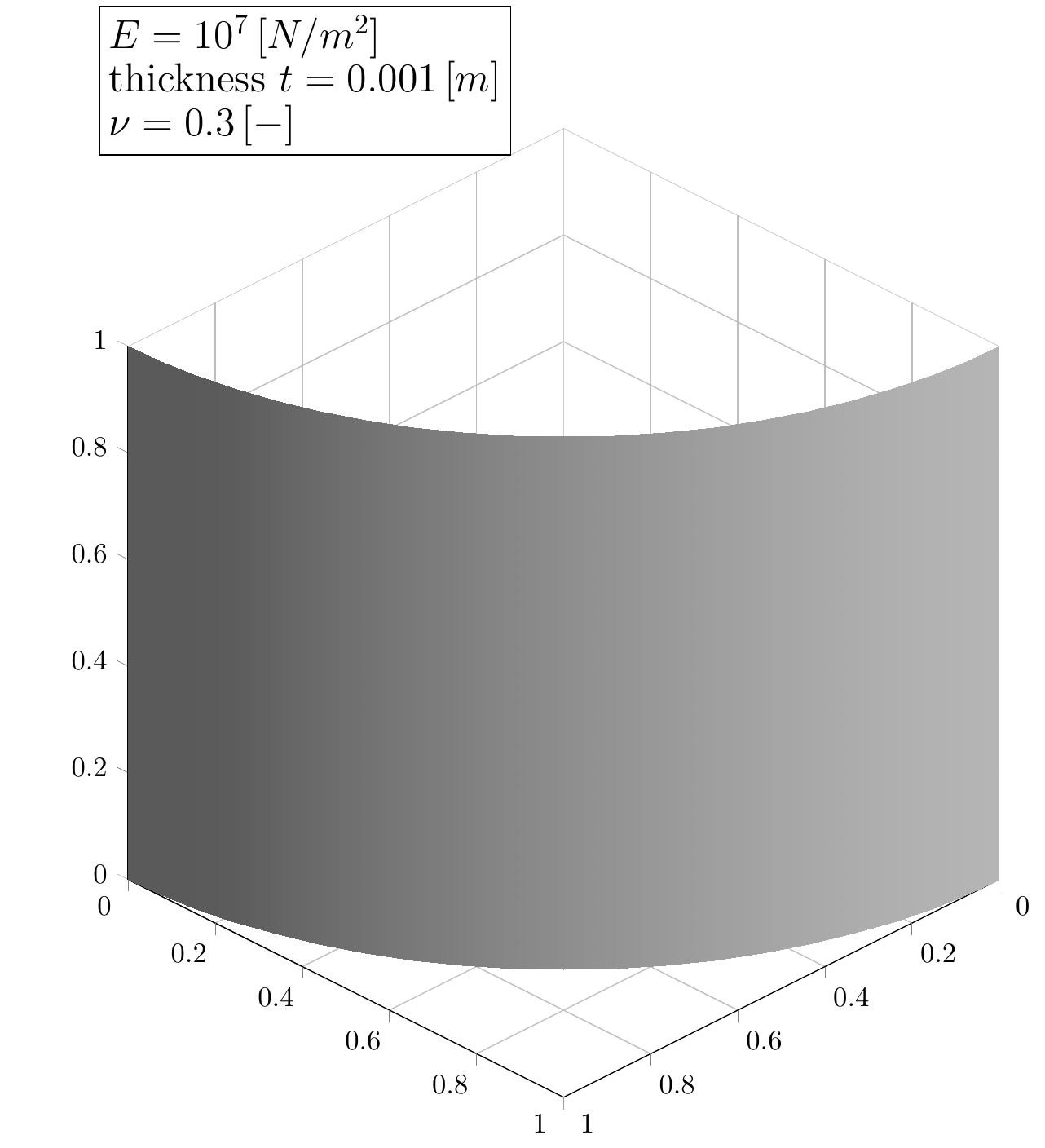}
		\caption{Geometry setup and physical parameters.}
	\end{subfigure}
	\hfill
	\begin{subfigure}[t]{0.495\textwidth}
		\centering
		\includegraphics[width=\textwidth]{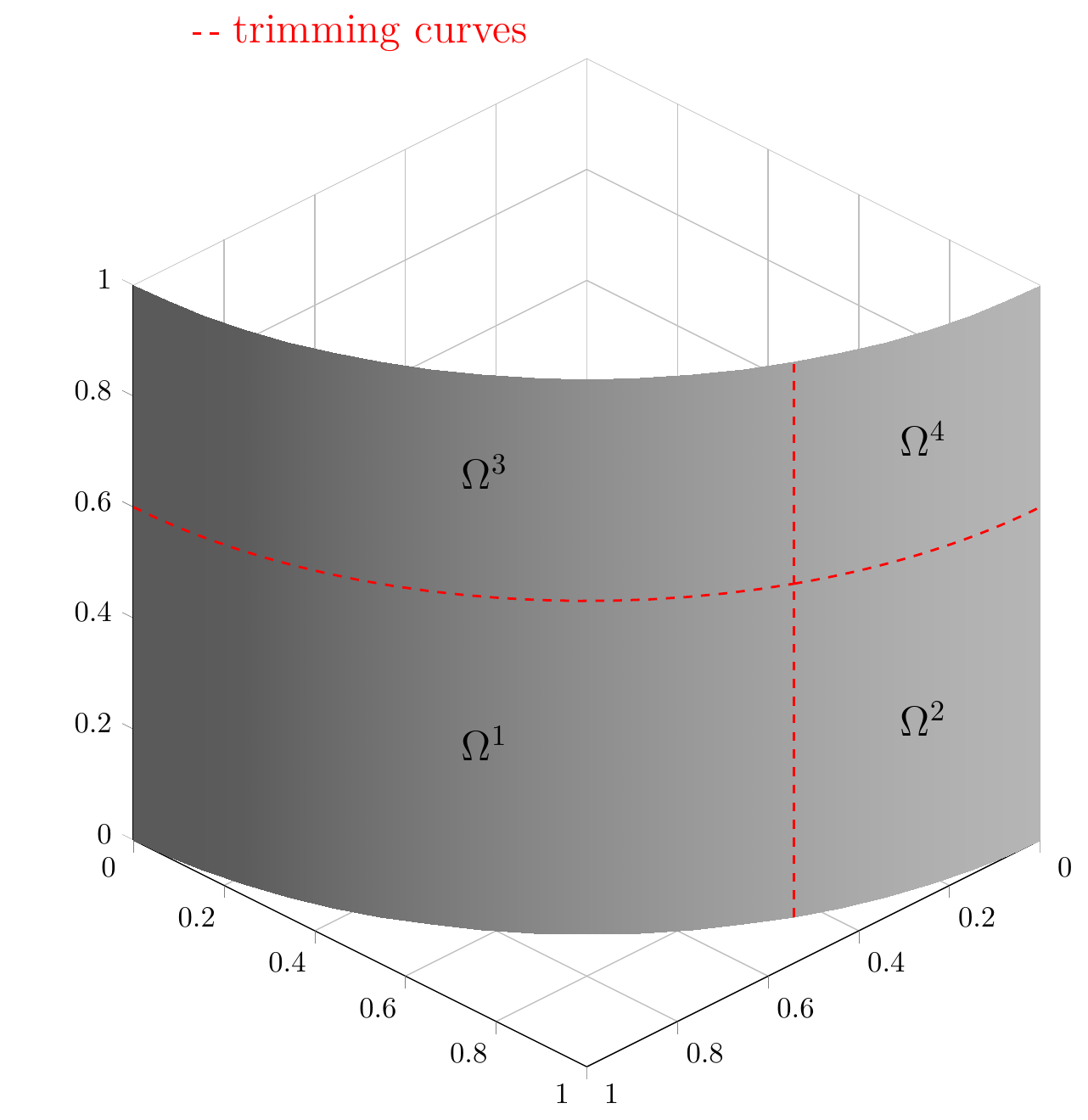}
		\caption{Trimmed multi-patch geometry.}
	\end{subfigure}	
	\caption{Problem setup and initial multi-patch geometry for the trimmed cylinder example.}
	\label{fig:setup_cylinder}
\end{figure}

\begin{figure}
	\centering
	\begin{subfigure}[t]{0.495\textwidth}
		\centering
		\input{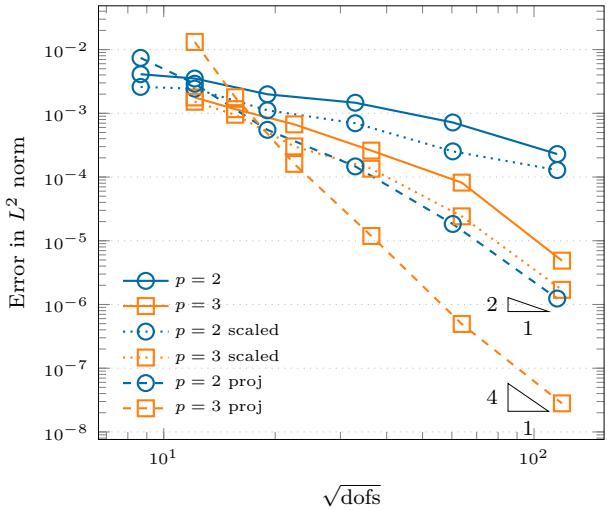}
		\caption{$L^2$ norm.}
	\end{subfigure}
	\hfill
	\begin{subfigure}[t]{0.495\textwidth}
		\centering
		\input{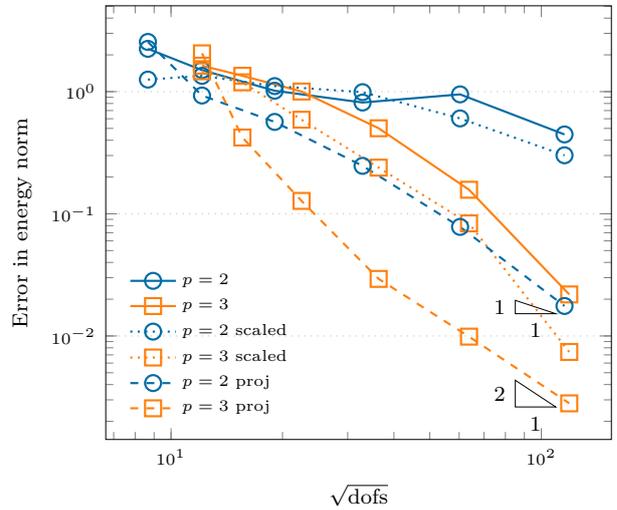}
		\caption{Energy norm.}
	\end{subfigure}	
	\caption{Convergence study of the error measured in the $L^2$ and energy norms for the trimmed cylinder example, B-splines of degree $p=2,3$. Comparison of a classic penalty method, the scaled version with respect to the problem parameters proposed in~\citep{Herrema2019} (\textit{scaled}) and our projection approach (\textit{proj}).}
	\label{fig:convergence_cylinder_t001}
\end{figure}

This example is again adapted from the shell obstacle course presented in~\citep{Benzaken2020}. We consider the computational domain $\Omega$ split into four trimmed subdomains $\Omega^i, i=1,\ldots,4$ as depicted in~\Cref{fig:setup_cylinder}. The corresponding trimmed interfaces $\gamma^\ell, \ell=1,\ldots,4$ are defined as quadratic B-spline curves. This numerical experiment tests the applicability of the proposed methodology to the coupling of trimmed multi-patch surfaces in the presence of cross-points.
Similarly to previous examples, the initial internal knots of patches $\Omega^2$ and $\Omega^3$ have been shifted by a factor $\sqrt{2}/100$ to achieve non-conforming discretization at the corresponding trimmed interfaces. Then, we set the Young's modulus, the Poisson's ratio and the thickness of the cylinder to $10^7 \, [Pa]$, $0.3 \, [-]$ and $0.001 \, [m]$, respectively.
The load function and boundary data are computed from the following manufactured solution:
\begin{align}
	\bm{u}^{\text{ex}} (\xi, \eta) = -(\xi - 1)^2 \xi^2 \eta (\eta - 1) \mathbf{a_3} \, ,
\end{align}
where $\mathbf{a_3}$ denotes the covariant vector in the thickness direction.
The convergence results for the error measured in the $L^2$ and energy norms, respectively, are depicted in~\Cref{fig:convergence_cylinder_t001}. Similarly to our previous findings, we observe a faster convergence behavior of the projection method in the pre-asymptotic regime, where interface locking is avoided on very coarse meshes. This results in a substantial gain of accuracy per-degree-of-freedom, which is particularly noticeable for quadratic B-splines.

\subsection{The DTU 10 MW Reference wind turbine blade}
\begin{figure}
	\centering
	\includegraphics[width=0.7\textwidth]{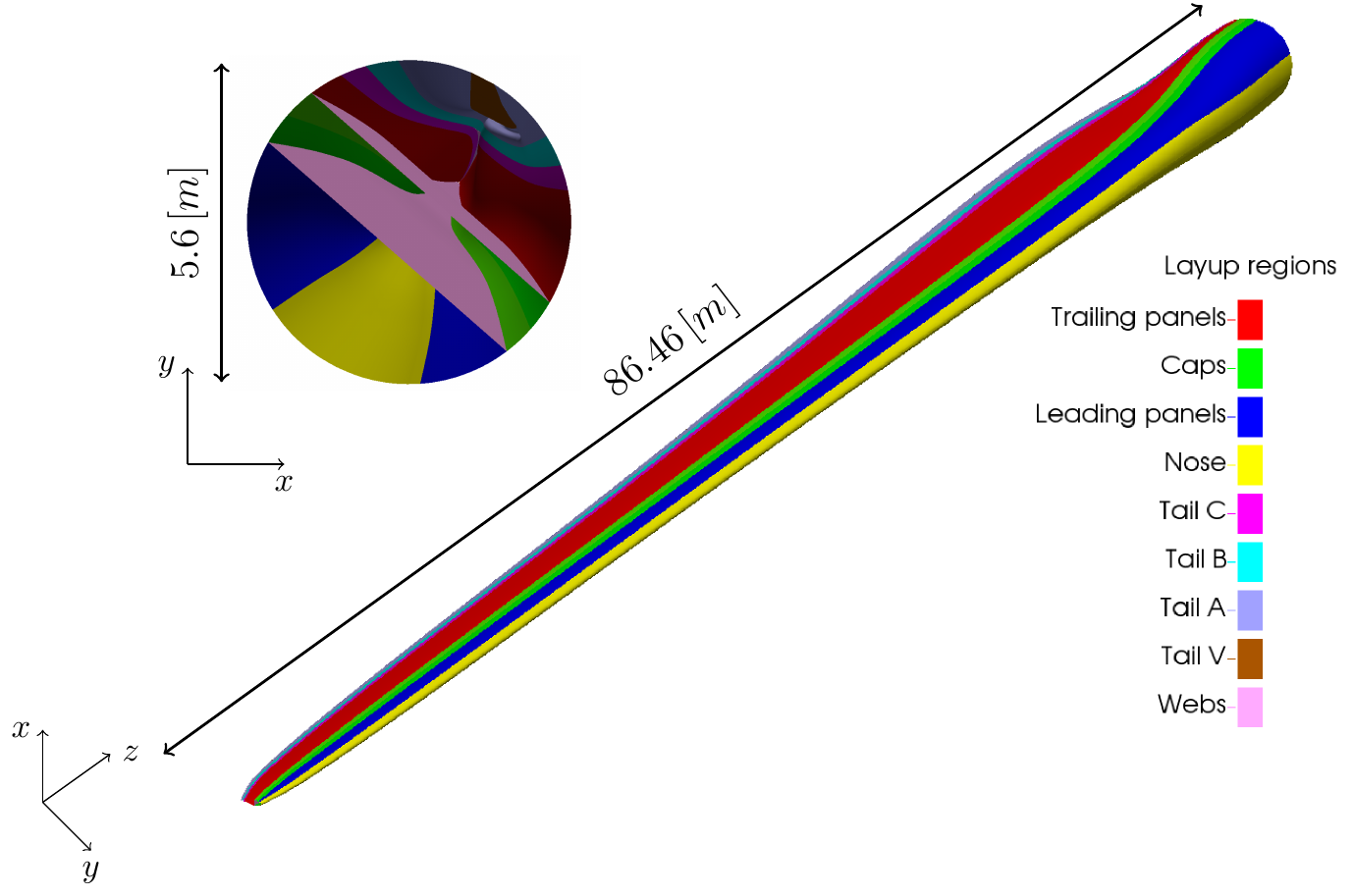}
	\caption{DTU 10 MW blade geometry and depiction of the circumferential regions used for the composite materials definition.}
	\label{fig:wind_turbine_regions}
\end{figure}

\begin{figure}
	\centering
	\includegraphics[width=0.9\textwidth]{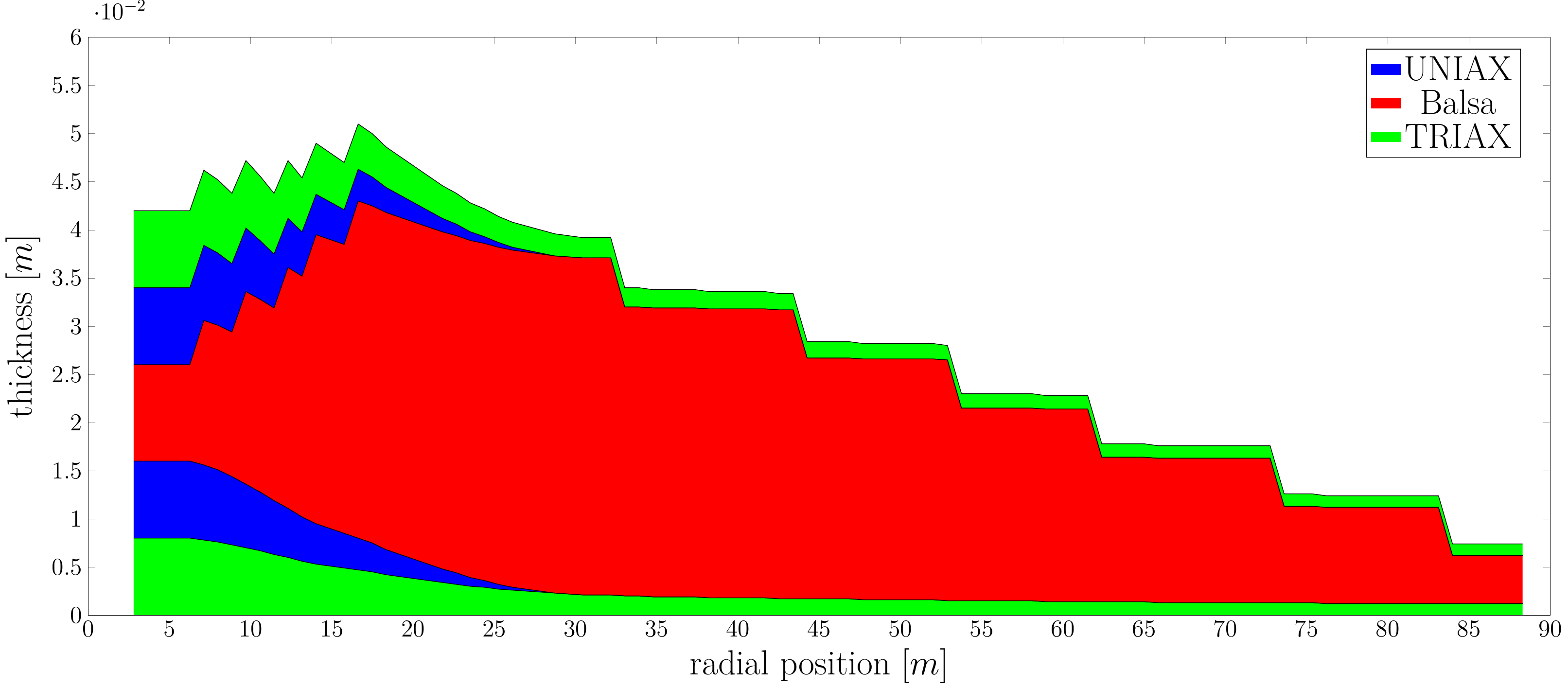}
	\caption{Composite layup along the radial direction of the leading panels.}
	\label{fig:wind_turbine_layup}
\end{figure}

\begin{figure}
	\centering
	\includegraphics[width=0.7\textwidth]{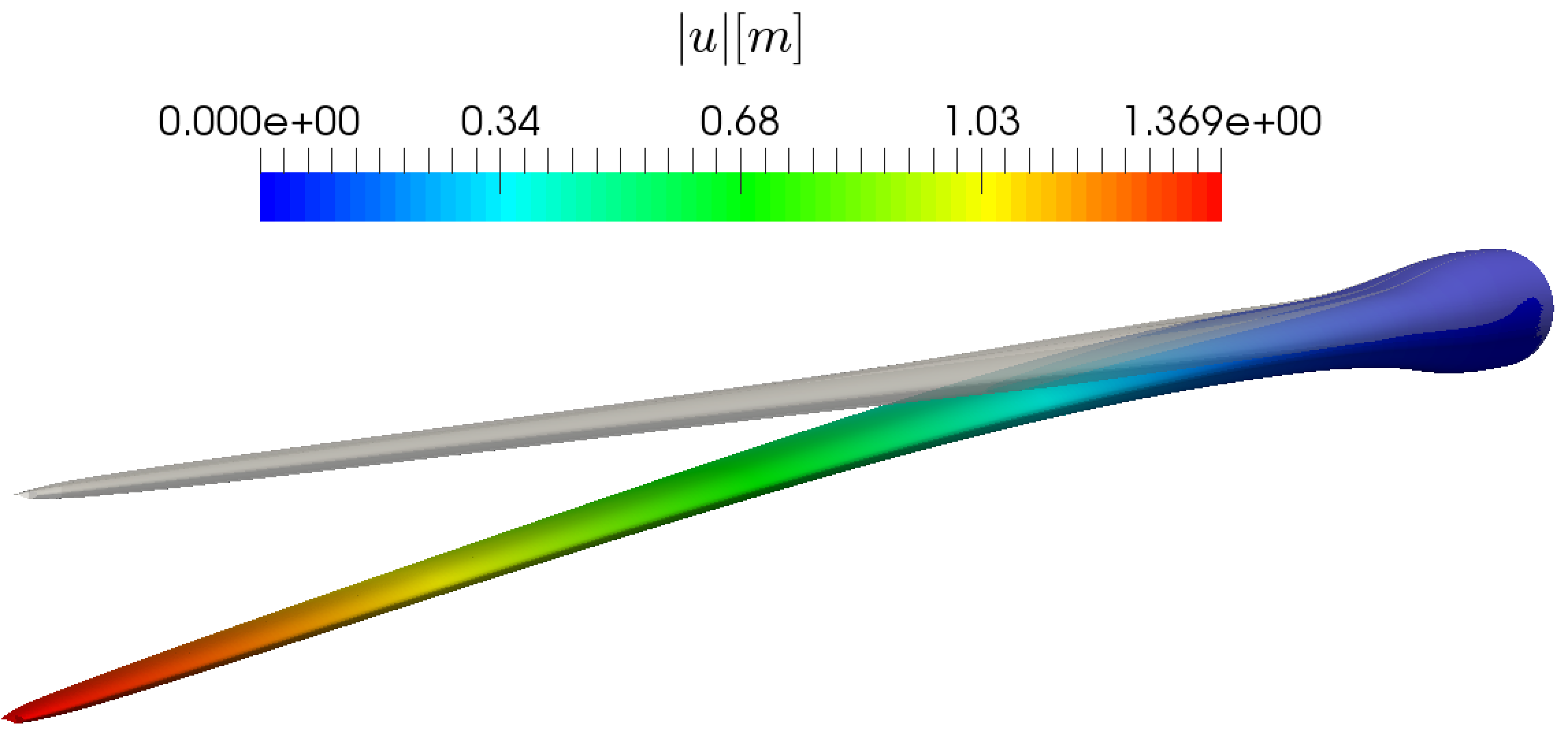}
	\caption{$L^2$ norm of the displacement field on the DTU 10 MW blade subject to gravitational load, B-splines of degree $p=2$, 89528 elements, results warped by a factor 5.}
	\label{fig:wind_turbine_disp}
\end{figure}

In our last example, we perform an isogeometric shell analysis of the DTU 10 MW Reference wind turbine blade~\citep{Bak2013}, whose design was inspired by the NREL 5 MW reference wind turbine~\citep{Jonkman2009}. The blade is modeled by 20 non-conforming cubic spline surfaces.
As noted in~\citep{Herrema2019}, a multi-patch design allows to accurately resolve material discontinuities along the patch interfaces. The outer shell of the blade and the internal shear webs are depicted in~\Cref{fig:wind_turbine_regions}. In the same figure, colored regions are used to define the corresponding composite layup, where each region has a different multi-directional ply stacking sequence and a varying thickness distribution along the spanwise direction. We summarize the most relevant mechanical properties in~\Cref{tab:mat_turbine}. Moreover, in~\Cref{fig:wind_turbine_layup}, we show the composite layup of the leading panels through the thickness as a function of the spanwise coordinate. For further details on material properties and thickness profiles we refer to~\citep{Bak2013}.
\begin{table}
	\centering
\begin{tabular}{lccccc} 
	Multi-directional ply & Uniax & Biax & Triax & Balsa &\\
	\hline Young's modulus $E_{1}$ & 41.63 & 13.92 & 21.79 & 0.050 & $[\mathrm{GPa}]$ \\
	Young's modulus $E_{2}$ & 14.93 & 13.92 & 14.67 & 0.050 & $[\mathrm{GPa}]$ \\
	Shear modulus $G_{12}$ & 5.047 & 11.50 & 9.413 & 0.01667 & $[\mathrm{GPa}]$ \\
	Poisson's ratio $\nu_{12}$ & 0.241 & 0.533 & 0.478 & 0.5 & $[-]$ \\
	Shear modulus $G_{13}=G_{23}$ & 5.04698 & 4.53864 & 4.53864 & 0.150 & $[\mathrm{GPa}]$ \\
	Mass density $\rho$ & 1915.5 & 1845.0 & 1845.0 & 110 & $[\mathrm{kg} / \mathrm{m}^{3}]$ \\
	\hline
\end{tabular}
\caption{Mechanical properties of the multi-directional plies.}\label{tab:mat_turbine}
\end{table}
For the analysis, we consider the response of the blade under gravity load, where the blade is modeled as clamped on the rotor side. The $L^2$ norm of the displacement field and the corresponding deflection of the blade are depicted in~\Cref{fig:wind_turbine_disp}. The results are obtained by employing a discretization of quadratic B-splines defined on 89528 elements.
Note that we directly import the geometry used in~\citep{Bak2013} for the structural analysis. 
\begin{remark}
	The latter is true for every patch except for webs A, B and C, which are obtained by linear extrusion of a generating spline, meaning that one linear element suffices to exactly describe the surface along the corresponding parametric direction. Therefore, $h$-refinement is performed along the latter direction by introducing 50 equidistributed knots. 
\end{remark} 
\noindent We remark that all the results on the blade have been obtained by setting the scaling factor $\beta = p$ in the penalty terms to limit the impact of the latter on the condition number of the stiffness matrix.

\subsubsection{Simplified topology optimization of webs A and B}
\begin{figure}
	\centering
	\includegraphics[width=\textwidth]{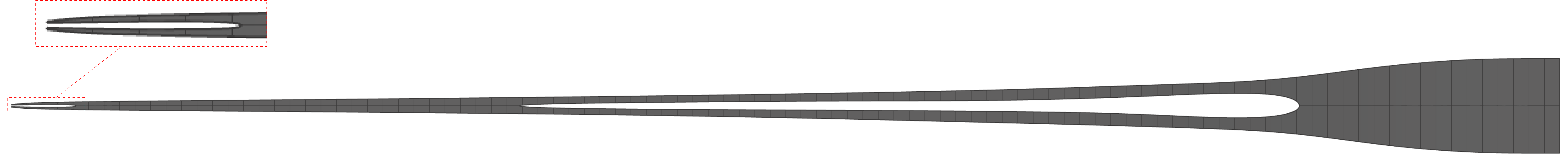}
	\caption{Design of web A after simplified topology optimization. Zoom on the detail at the end of the blade.}
	\label{fig:web_opt}
\end{figure}

\begin{figure}
	\centering
	\begin{subfigure}[t]{0.495\textwidth}
	\centering
	\includegraphics[width=\textwidth]{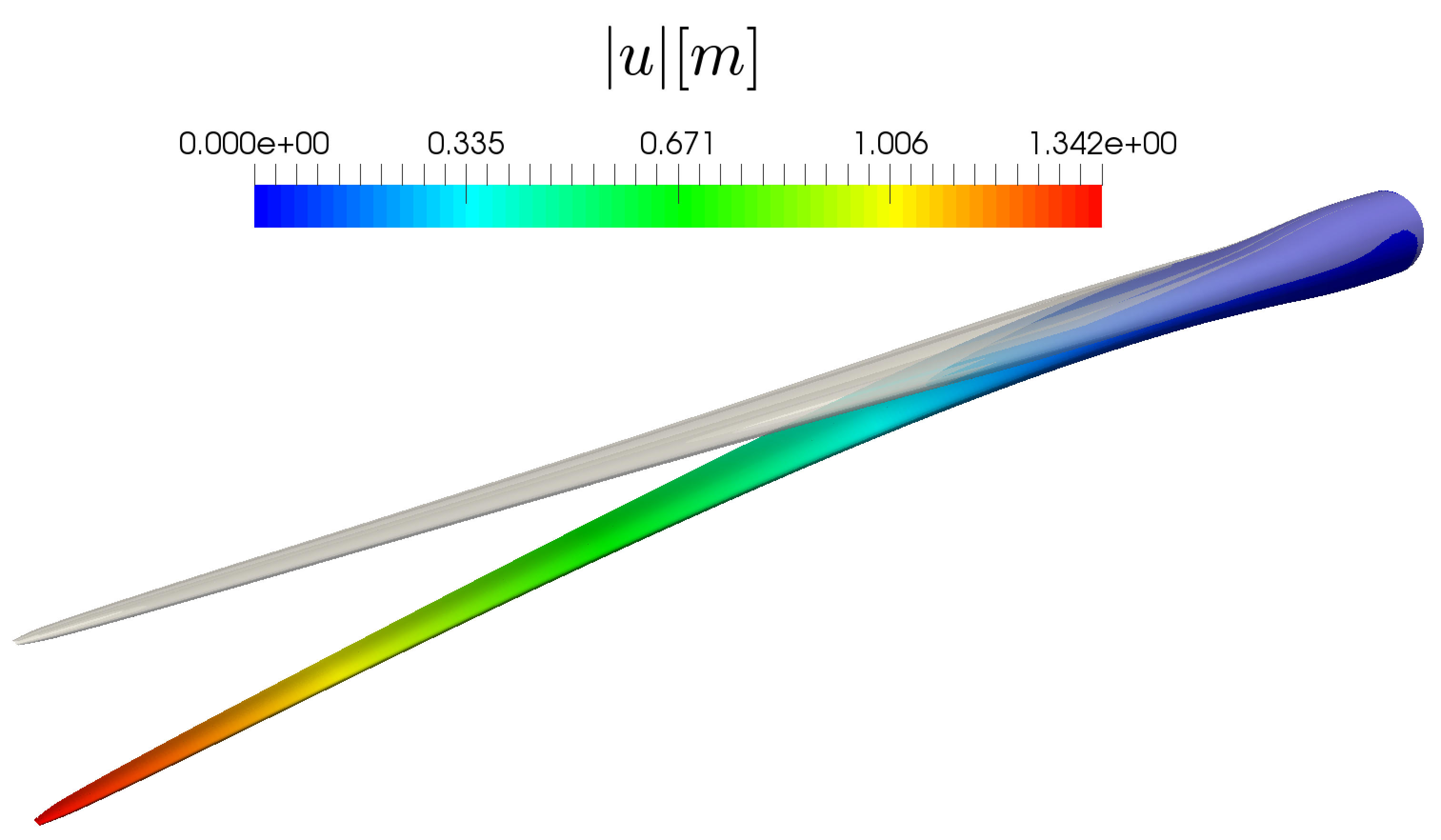}
	\caption{Full blade.}
	\end{subfigure}
	\hfill
		\begin{subfigure}[t]{0.495\textwidth}
		\centering
		\includegraphics[width=\textwidth]{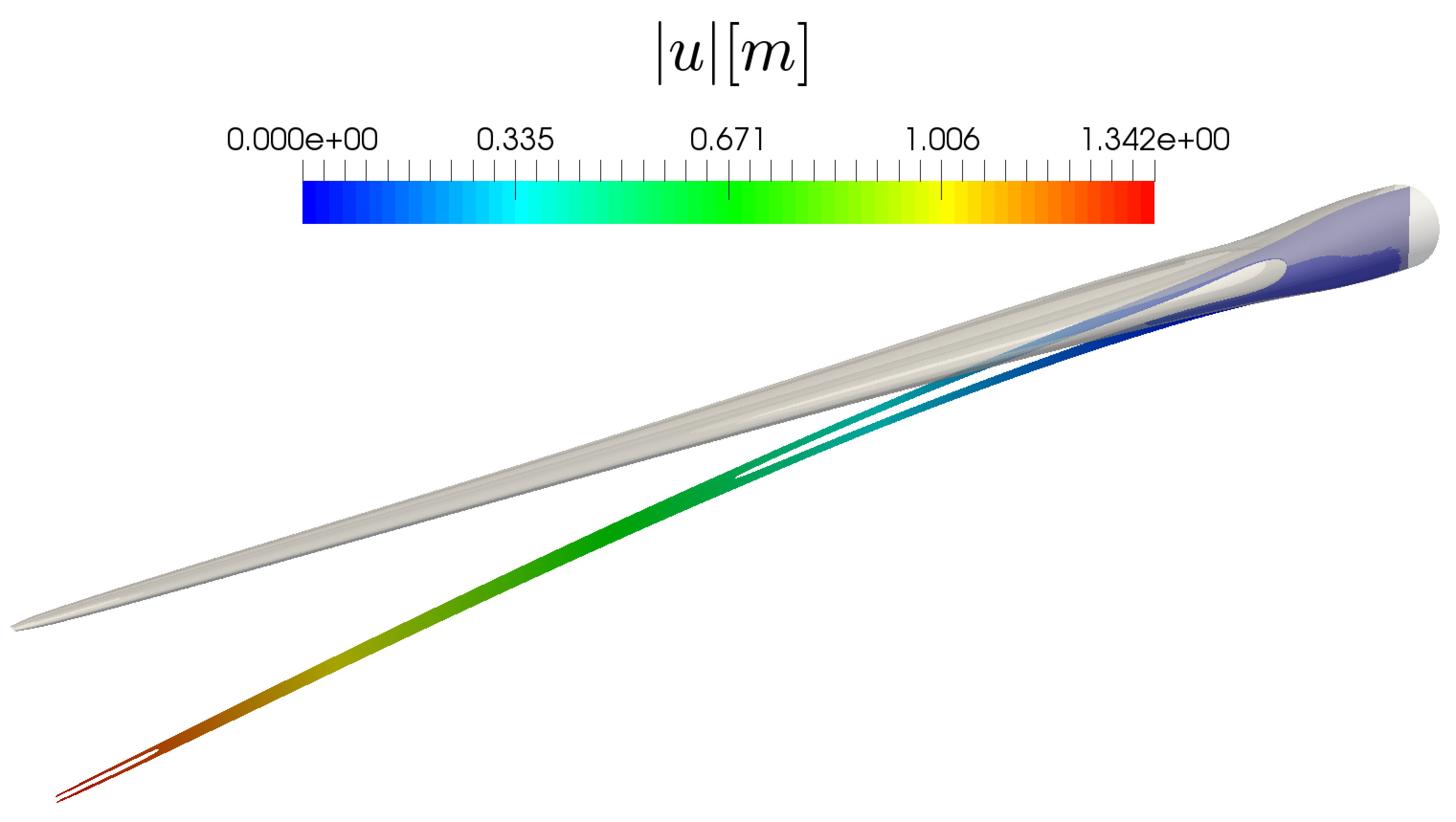}
		\caption{Web A.}
	\end{subfigure}
	\caption{$L^2$ norm of the displacement field on the topologically optimized blade subject to gravitational load, B-splines of degree $p=2$, 84250 elements, results warped by a factor 5.}
	\label{fig:wind_turbine_opt_disp}
\end{figure}

This example is meant to show the applicability of the proposed methodology to trimmed geometries obtained by a simplified topology optimization. Note that this numerical test is just a showcase of the flexibility of our computational framework and a realistic topology optimization of the webs is beyond the scope of this work. Furthermore, the geometric operation described in this section are based on an heuristic engineering approach.
We trimmed away from the original geometry two holes, one close to the center of the structure and another at the end the web. This design is obtained by adapting the optimized solution presented in~\citep{Albanesi2020}, where the final geometry is depicted in~\Cref{fig:web_opt}. This design results in a reduction of $\approx 20.6 \%$ of the original total mass of web A. Similarly, we perform the same operations on web B.
The $L^2$ norm of the displacement field and the corresponding deflection of the blade with trimmed webs are depicted in~\Cref{fig:wind_turbine_opt_disp}, where we observe a reduction of $\approx 2.0 \%$ in tip displacement related to the loss of total mass of the structure. The results are obtained with quadratic B-splines defined on a total of 84250 active elements.

\section{Conclusions} \label{sec:conclusions}

In this contribution we have extended the methodology presented in~\citep{Coradello2020projected} to the coupling of trimmed, non-matching surfaces in the context of isogeometric Kirchhoff-Love shells. The strategy is based on the $L^2$-projection of suitable penalty terms at the corresponding coupling interface onto a reduced space of degree $p-2$ with respect to the chosen approximation space. On one hand, the projection mitigates the detrimental effects related to interface locking starting from very coarse discretization. On the other hand, it gives us insights into the proper scaling of the penalty parameters based on the underlying discretization. Consequently, the proposed coupling method retains the optimal rates of convergence achievable by B-splines, as demonstrated by our findings on an extensive series of benchmark problems. Moreover, building upon the penalty factors studied in~\citep{Herrema2019}, our approach is fully parameter-free, since the penalty coefficients are completely determined by the problem setup. Similarly to what was observed in~\citep{Coradello2020projected}, our strategy is particularly suited for spline spaces of moderate degrees $p=2,3$, where the projection turns out to be computationally efficient and the condition number stemming from the super-penalty does not yield a significant deterioration of the solution accuracy. 
Then, the applicability of our method to tackle complex, industrial-like structure has been studied. In particular, we have performed a static shell analysis of the DTU 10MW Reference wind turbine blade~\citep{Bak2013}, whose design is composed by 20 non-conforming cubic patches. Furthermore, since trimming is naturally incorporated into our methodology, we have carried out a simplified topology optimization of the internal webs. This example demonstrates that the proposed computational framework is readily applicable to tackle industrial optimization loops of engineering relevance.  

To conclude, we have numerically validated the wide range of applicability and robustness of the proposed projected super-penalty approach for coupling trimmed Kirchhoff-Love shells, where the method does not suffer from boundary locking and the penalty parameters are automatically derived from the problem setup to attain optimal rates of convergence.
Finally, in this contribution only a geometrically linear shell formulation and its static behavior have been considered. Potential future research directions include the extension of the proposed approach to account for geometrical and material non-linearities and the assessment of its performance in dynamics.

\section*{Acknowledgements} 
The authors L.~Coradello and A.~Buffa gratefully acknowledge the support of the European Research Council, via the ERC AdG project CHANGE n.694515. The author A.~Buffa also gratefully acknowledges the support of the H2020-FetOpen-Ria project ADAM$^2$ n.862025. The author J.~Kiendl gratefully acknowledges the support of the European Research Council through the ERC Consolidator Grant FDM$^2$, grant n.864482.

\bibliographystyle{plainnatnourl}
\bibliography{library.bib}

\end{document}